\documentclass[10pt]{article}

 \usepackage[utf8]{inputenc}
\usepackage[english]{babel}
\usepackage{amsmath, amsfonts, amstext, amssymb} 

 \usepackage{makeidx,epsfig,lscape}
 \usepackage{color,colortbl}
 \usepackage{fancyhdr}
 \usepackage{xcolor,pict2e}

\usepackage{graphicx}

\newcommand{\N}{\mathbb{N}}
\newcommand{\Z}{\mathbb{Z}}
\newcommand{\R}{\mathbb{R}}

\newcommand{\A}{\mathcal{A}}
\newcommand{\la}{L_{\beta}}
\newcommand{\C}{\mathfrak{C}}

\newcommand{\Di}{\Delta^{(i)}}

\newcommand{\s}{\tilde{S}_{-\beta}}
\newtheorem{theo}{Theorem}
\newtheorem{coro}{Corollary}
\newtheorem{lem}{Lemma}
\newtheorem{prop}{Proposition}

\newtheorem{rem}{Remark}

\newenvironment{preu}{\textbf{Proof }\rm}{\par\hfill$\square$}
\newenvironment{preu 1}{\textbf{Proof of Theorem \ref{t2}.}\rm}{\par\hfill$\square$}
\newenvironment{preu 2}{\textbf{Proof of Theorem \ref{t1}.}\rm}{\par\hfill$\square$}
\newenvironment{preu 3}{\textbf{Proof of Theorem \ref{th3}}\rm}{\par\hfill$\square$}
\newenvironment{preu 4}{\textbf{Proof of Theorem \ref{theo3}}\rm}{\par\hfill$\square$}
\newenvironment{preu 5}{\textbf{Proof of Theorem \ref{p2}}\rm}{\par\hfill$\square$}
\newtheorem{defi}{Definition}
\newtheorem{exple}{Example}

\usepackage{amssymb}
\usepackage{datetime}
\usepackage{lineno}

 \renewcommand{\headrulewidth}{0pt}
 \renewcommand{\footrulewidth}{0.5pt}

 \definecolor{myaqua}{rgb}{0.0,0.5,0.55}
 \definecolor{lightaqua}{rgb}{0.75,0.95,0.95}

 \usepackage[colorlinks = true,
            linkcolor = myaqua,
            urlcolor  = blue,
            citecolor = red,
            pdfpagemode=UseOutlines,
            bookmarksnumbered=true,
            pdfpagelabels,
            breaklinks]{hyperref}

\usepackage{caption}
\usepackage{floatrow}

\def\lin#1#2{\textcolor[rgb]{0.6,0.6,0.6}{\vspace*{#1mm} \hrule
   height 3 pt \vspace*{#2mm}}}
\def\bt{\begin{tabular}}
\def\et{\end{tabular}}
\def\and{\mbox{ and }}

\def\1{{\bf 1}}

 \def\boxx#1#2#3#4#5{
 {\linethickness{#4pt}\put(#1,#5){\color{myaqua}{\line(1,0){#3}}}}
 \multiput(#1,#2)(0,#4){2}{\line(1,0){#3}}
 \multiput(#1,#2)(#3,0){2}{\line(0,1){#4}}
  }

\begin{document}


 $\mbox{ }$

 \vskip 12mm
 
 { 

{\noindent{\Large\bf\color{myaqua}
  The $(-\beta)$-shift and associated Zeta Function  }} 
%
\\[6mm]
{\bf Florent NGUEMA NDONG }}
\\[2mm]
{ 
 $^1$ Universit\'e des Sciences et Techniques de Masuku 
 \\
Email: \href{mailto:florentnn@yahoo.fr}{\color{blue}{\underline{\smash{florentnn@yahoo.fr}}}}\\[1mm]
\lin{5}{7}

 {  
 {\noindent{\large\bf\color{myaqua} Abstract}{\bf \\[3mm]
 \textup{
    Given a real number $ \beta > 1$, we study the associated $ (-\beta)$-shift introduced by S. Ito and T. Sadahiro. We compare some aspects of the $(-\beta)$-shift to the $\beta$-shift. When the expansion in base $ -\beta $ of $ -\frac{\beta}{\beta+1} $ is periodic with odd period or when $ \beta $ is strictly less than the golden ratio, the $ (-\beta)$-shift cannot be coded because its language is not transitive. This intransitivity of words explains the existence of gaps in the interval $[-\frac{\beta}{\beta+1}, \frac{1}{\beta+1})$. We observe that an intransitive word appears in the $(-\beta)$-expansion of a real number taken in the gap. Furthermore, we determine the Zeta function $\zeta_{-\beta}$ of the $(-\beta)$-transformation and the associated lap-counting function $L_{T_{-\beta}}$. These two functions are related by $ \zeta_{-\beta}=(1-z^2)L_{T_{-\beta}}$. We observe some similarities with the zeta function of the $\beta$-transformation. The function $\zeta_{-\beta}$ is meromorphic in the unit disk, is holomorphic in the open disk $ \{z \hspace{0.1cm}: \hspace{0.1cm} \vert z \vert < \frac{1}{\beta} \}$, has a simple pole at $ \frac{1}{\beta}$ and no other singularities $ z $ such that $\vert z \vert = \frac{1}{\beta}$. We also note an influence of gaps ($\beta$ less than the golden ratio) on the zeta function. In  factors of the denominator of $\zeta_{-\beta}$, the coefficients count the words generating gaps.
 }}}
 \\[4mm]
 {\noindent{\large\bf\color{myaqua} Keywords}{\bf \\[3mm]
  Negative basis; $\beta$-expansions; coded system; transitive system; zeta function. 
}}}
\lin{3}{1}

\renewcommand{\headrulewidth}{0.5pt}
\renewcommand{\footrulewidth}{0pt}

 \pagestyle{fancy}
 \fancyfoot{}
 \fancyhead{} 
 \fancyhf{}
 \fancyhead[RO]{\leavevmode \put(-90,0){\color{myaqua}F. NGUEMA NDONG} \boxx{15}{-10}{10}{50}{15} }
 \fancyfoot[C]{\leavevmode
 \put(-2.5,-3){\color{myaqua}\thepage}}

 \renewcommand{\headrule}{\hbox to\headwidth{\color{myaqua}\leaders\hrule height \headrulewidth\hfill}}

 \section{Introduction}

{ 
\selectfont
 \noindent 
%

The $\beta$-transformation has been extensively studied since the seminal paper of R\'enyi in 1957 (see \cite{MR0097374}).
There is a huge amount of literature on the map itself and on the associated symbolic dynamics.
Over the past decade, people became interested in the $(-\beta)$-transformation, changing the slope of the map from positive to negative. Various studies have focused on the similarities and differences
between the two maps from several points of view. This paper fits in this last line of research.

The paper compares two aspects of the $(-\beta)$-shift to the $\beta$-shift. For the $\beta$-shift it is known that a prefix code always exists. The paper first investigates whether or not the same is true for the $(-\beta)$-shift.  By $ (-\beta)$-shift (respectively $\beta$-shift)  we understand the closure of the set of expansions in base $-\beta $ (respectively $ \beta$). The conclusion is stated in Theorem \ref{t1}, which says that a prefix code exists in the negative case only under certain conditions, namely if and only if $\beta$ is bigger than the golden ratio and the orbit of the left endpoint of the domain of the $(-\beta)$-transformation is not periodic with odd period. It turns out that the discovered prefix codes are recurrent. Even though the codes can be complicated, the zeta functions apparently have a very simple form (see Theorem \ref{t2}) and it can be noted the similarities with that of the positive parameter determined in \cite{MR0176482}, \cite{MR1279470} and \cite{MR1710789}.



R\'enyi introduced the $ \beta$-expansion of positive real numbers in \cite{MR0097374}: for a fixed real $ \beta > 1 $, all non-negative real number $ x $ has one representation in base $ \beta $. He defined the $ \beta$-transformation $ T_{\beta} $ from $ [0, 1) $ into $ [0, 1) $ by
\begin{equation*}
  T_\beta (y) = \beta y -\lfloor \beta y\rfloor 
 \end{equation*}
 where $ \lfloor x \rfloor $ denotes the largest integer less than $ x $. We can find a sequence of positive integers $ ( x_i)_{i \geq -n+1} $ such that
 $ x = \underset{i \geq -n+1}{\sum}\dfrac{x_i}{\beta^i} $ where $x_i \in  \left\lbrace 0, 1, \cdots, \lfloor \beta \rfloor \right\rbrace $, $ x_{-n+i} = \lfloor \beta T^{i-1}_\beta(\frac{x}{\beta^n}) \rfloor $ and $ n $ is the smallest positive integer for which $ \frac{x}{\beta^n} $ belongs to the interval $ [0, 1) $. Various topics about $ \beta$-expansion have been studied.

Let $ b $ be an integer strictly bigger than $ 1 $. It is well-known that any number can be represented without a sign in base $ -b $ thanks to digits 
of the alphabet $ \left\lbrace  0, 1, \cdots, b-1 \right\rbrace $.
 In 2009, S. Ito and T. Sadahiro in \cite{MR2534912} generalized this approach for all negative base (integer or not).
 They defined a $(-\beta)$-transformation $ T_{-\beta} $, the map from the interval $ I_{\beta} = [-\frac{\beta}{\beta+1}, \frac{1}{\beta}) $ into itself such that
\begin{equation*}
  T_{-\beta}(x) = -\beta x - \lfloor -\beta x + \frac{\beta}{\beta+1} \rfloor.
\end{equation*}
The expansion in base $ -\beta $ of a real $ x $ (denoted by $ d(x, -\beta )$) is given by the following algorithm: 
\begin{itemize}
\item if  $ x $ belongs to $ I_\beta $, $ d(x, -\beta) = . x_1x_2\cdots, $ where  
\begin{equation*}
  x_i = \lfloor -\beta T_{-\beta}^{i-1}(x) + \frac{\beta}{\beta+1} \rfloor; 
 \end{equation*}
\item if $ x \not\in I_{\beta} $, one finds the smallest integer $ n $ 
for which one has $ \frac{x}{(-\beta)^n} \in I_{\beta}$. In this case, the expansion is 
$ d(x, -\beta) =x_{-n+1}\cdots x_0. x_1x_2 \cdots $, where
\begin{equation*}
  x_{-n+i} = \lfloor -\beta T_{-\beta}^{i-1}\left(\frac{x}{(-\beta)^n}\right)+\frac{\beta}{\beta+1} \rfloor,\text{ $ i \geq 1 $ }.
\end{equation*}
\end{itemize}
In both cases, $ x = \underset{ i \geq -n+1}{\sum}\dfrac{x_i}{(-\beta)^i} $. 
If there is no mixup, we often denote $ d(x, -\beta)$ by $ (x_i)_{i \geq -n+1} $.

\subsection*{Alternating lexicographic order}
\vspace*{0.3cm}
\begin{defi}
Let $\A = \{0, 1, \cdots, k \}$ be an alphabet. A word on $\A$ is a concatenation of elements of $\A$. Let $ x=x_1 x_2 \cdots x_n $ be a word on $ \A $, we call sub-word (or factor) of $x$ all word which appears in $x$. 
$\A^*$ denotes the set of words on $ \A $.
\end{defi}
Let $ \A = \{0, 1, \cdots, d_1 \} $. We endow $ \A^*$ with the order $\prec $ defined as follows: for a pair of words  $ X_n= x_1 x_2 \cdots x_n $ and $ Y_n = y_1 y_2 \cdots y_n $, $n \in \N^{*}$
\begin{equation*}
X_n \prec Y_n \Leftrightarrow \exists k, \text{ $ x_1 \cdots x_{k-1} = y_1 \cdots y_{k-1} $ and $ (-1)^k(x_k-y_k)< 0 $}.
\end{equation*}
We say $ X_n $ is less than $ Y_n $ with respect to $ \prec $. Comparison is also possible for words of distinct lengths (see \cite{NguemaNdong20161}).

The order $``\prec"$ is called \textit{alternating lexicographic order}. With this order, we have the following property: let $ u $, $ v $ and $ w $ be three words on $ \A $,
\begin{equation}
u \prec v \Rightarrow \begin{cases}
                       w u \prec w v &\text{ if $ |w| $ is even }\\
											 w v \prec w u &\text{ if $ |w| $ is odd }, \\
											 u w \prec v w &\text{ in all cases}. 
											\end{cases} \label{ef1}
\end{equation}

\subsection*{The $(-\beta)$-shift}
The \textit{$ (-\beta)$-shift } is the closure of the set of expansions in base $-\beta$. 
The sequence $d(-\frac{\beta}{\beta+1}, -\beta)$ plays an important role in the characterization of this set.
 In the following, 
 \begin{equation*}
  d(l_\beta, -\beta) = \cdot d_1 d_2 \cdots 
 \end{equation*}
  with  $ l_\beta = -\frac{\beta}{\beta + 1} $, $ r_\beta =  \frac{1}{\beta + 1} $, $ d_0 = 0 $ and
 \begin{equation}
 (d_{i}^{*})_{i \geq 1} = \begin{cases}
                         \overline{(d_1, \cdots, d_{2n_i-2}, d_{2n-1}-1, 0)} &\text{ if $ (d_i)_{i \geq 1} = \overline{(d_1, \cdots, d_{2n-1} )} $ }\\
                         (d_i)_{i \geq 1} &\text{ otherwise } 
                         \end{cases}\label{D}
\end{equation}      
where $ \overline{t} = ttt \cdots $ is the purely periodic sequence with period $ t $.
\begin{equation}
S_{-\beta} = \{ x_k x_{k+1} \cdots x_0 \cdot x_1\cdots \vert (d_i)_{i\geq 1} \preceq (x_i)_{i \geq m} \preceq (d_{i-1}^{*})_{i \geq 1}, \forall m\geq k, \forall k \} \label{S}.
\end{equation}

\section{Coded negative beta-shift}

Let us start by giving the definitions of the main terms used throughout this paper.

\subsection{Definitions}

\begin{defi}
Let $\A = \{0, 1, \cdots, k \}$ be an alphabet. 
A language $L$ on $\A$ is a set of words on $\A$ (or the set of finite sequences on $\A$).

A language $ L $ is extendable if for all word $ x_1x_2 \cdots x_n $ in $ L $, there exist two letters $ a $ and $ b $ in $ \A $ such that $ ax_1 x_2\cdots x_n b $ belongs to $ L $. It is said transitive if $ \forall v, w \in L $, there exists $ u $ such that $ vuw \in L $.
\end{defi}

Let $ \A = \{0, 1, \cdots, k \}$ be a finite alphabet. One endows $\A^{\N}$ (or $ \A^{\Z}$) with the topology product of the discrete topology on $\A$. Let $\sigma $ be the map from $ \A^{\N} $ (or $\A^{\Z}$) into itself defined by $\sigma((x_n)_{n \in \N}) = (x_{n+1})_{n \in \N}$. The closed $\sigma$-invariant subsets of $\A^{\N}$ are called sub-shifts. 

\begin{defi}
Let $S$ be a sub-shift on $\A$. The pair $ (S, \sigma)$ is called symbolic dynamical system. $ (S, \sigma) $ will be denoted by $ S $ whenever there is no risk of mixup.
\end{defi}

A language $ L_S $ of a dynamical system $ S $ is the set of factors of words of $ S $. The system is transitive if its language is transitive. 

\begin{defi}
 A code $Y$ on $ \A $ is a language such that, for any equality 
\[ x_1\cdots x_n = y_1 \cdots y_p \] 
with $ x_i, y_i \in Y $, one has $ n=p $ and $ x_i = y_i $.
\end{defi}

Let $ L$  be a language on $\A$. In the following, $L^{*}$ denotes the concatenations of words of $ L$.

\begin{defi}
 A prefix code is a language $ \C $ for which no word is the beginning of another. 
\begin{equation*}
  \forall x, y \in \C, \text{ $  x = yz \Rightarrow x = y \text{ and } z = \varepsilon $ }
\end{equation*}
 where $ \varepsilon $ is the empty word. 
 
 If in $ \C $, no word ends by another, then the language $ \C $ is a suffix code.

The symbolic dynamical  system $ S $ is said coded by the prefix code $\C$ if $L_S $ is the set of factors of words contained in $\C^{*}$. 
\end{defi}

\begin{defi} 

Let $ L $ be a language on $ \A $. The radius $ \rho_L$ of the power series $ \sum\limits_{ n \geq 1} card (L\cap \A^n)z^n $ is called radius of convergence of $ L $.

A prefix code $ \C $ is said recurrent positive if 
\begin{equation*}
\sum\limits_{x \in \C} \rho_{\C^{*}}^{ \vert x \vert } = 1 \text{ and } \sum\limits_{x \in \C} \vert x \vert \rho_{\C^{*}}^{\vert x \vert} < \infty.
\end{equation*}
\end{defi}
 
\subsection{Constructing of a code} 

Let $ \beta $ be a real number bigger than 1 and $ S_\beta $ the associated $\beta$-shift. If $ (a_i)_{i \geq 1} $ denotes the expansion of 1 in base $ \beta $, the $\beta$-shift, $ S_{\beta}$ is coded by the prefix code $ Y_{\beta} $ defined as follows: 
\begin{equation}
Y_{\beta} = \left\lbrace a_1a_2\cdots a_n j ; j < a_{n+1}, \forall n \in \N \right\rbrace.\label{BS}
\end{equation}
 Hence, all $ \beta $-shift is coded. It is one of the differences between $ \beta$-shifts and $ (-\beta)$-shifts. In fact, the $ (-\beta)$-shifts, given in \eqref{S} are not always coded. One of the natural and important question is whether the $(-\beta)$-shift is coded. In this section, we shall provide some contributions to this question. Furthermore, when it is coded, it is not easy to find a set of words coding its language. 
We distinguish two cases: for all $ i $, $0\leq d_{2i} < d_1 $ and $ d_{2i} = d_1 $ for some $ i $. 

Throughout in the rest of this paper, $ L_{\beta}$ denotes the language of the $(-\beta)$-shift.

The following theorem states the conditions on the parameter $ \beta $ to have a coded $(-\beta)$-shift.  
\begin{theo}\label{t1}
Let $ \beta $ be a real number greater than 1. The associated $(-\beta)$-shift $ S_{-\beta}$ is coded if only if $ \beta $ is greater than or equal to the golden ratio and $ d(l_\beta, -\beta)$ is not periodic with odd period.
\end{theo}
\begin{lem}\label{lem1}
Let $ \beta> 1$ and $ (d_i)_{i \geq 1}$ the $ (-\beta)$-expansion of $ -\frac{\beta}{\beta+1} $. If for all integer $ i $, $ d_{2i} < d_1 $, then $ \beta $ is bigger than or equal to the golden ratio.
\end{lem}

\begin{preu}
Let $ \gamma_0$ be the golden ratio. To prove this lemma, it is enough to determine  the $(-\gamma_0)$-expansion of $ l_{\gamma_0}=-\frac{\gamma_0}{1+\gamma_0}$. In fact, if we consider two real numbers $ \beta $ and $ \alpha $ strictly greater than 1 then,
\begin{equation*}
 d(l_\alpha, -\alpha) \prec d(l_\beta, -\beta) \Leftrightarrow \beta < \alpha.
\end{equation*}
We obtain the equivalence above thanks to Lemma 3 of \cite{NguemaNdong20161}. 
\begin{equation*}
d(l_{\gamma_0}, -\gamma_0) = \cdot 1 \overline{0}.
\end{equation*}

We assume that $ \beta < \gamma_0$ and $ d(l_\beta, -\beta) = (d_i)_{i \geq 1}$. Then, 
\begin{equation*}
1\overline{0} \prec d_1 d_2 d_3 \cdots. 
\end{equation*}
That means there exists $ n $ such that $ -(-1)^nd_n < 0 $ and $ d_1 d_2 \cdots d_{n-1} = 1 0\cdots 0 $. Thus, $ d_n = 1 $ and $ n $ is even. So, there exists $ i_0$  such that $n=2i_0$ and $ d_{2i_0} = 1= d_1 $. Thus, if $ d_{2i}<d_1 $ for all $ i $, $ \beta \geq \gamma_0 $. 
\end{preu}
\begin{prop}\label{prop1}
Let $ \beta $ be a real number greater than 1. We denote by $ S_{-\beta}$ the associated $(-\beta)$-shift, $ (d_i)_{i \geq 1}$ the $ (-\beta)$-expansion of $l_\beta=-\frac{\beta}{\beta+1} $. If $ (d_i)_{i \geq 1}$ is periodic with odd period or $ \beta < \frac{1+\sqrt{5}}{2} $ then, $ S_{-\beta}$ is not transitive. 
\end{prop}

\begin{preu}
\begin{itemize}
\item Assume $ \beta < \frac{1+\sqrt{5}}{2} $. By Lemma \ref{lem1}, there exists $ i_0$ such that 
\begin{equation*}
 d_1\cdots d_{2i_0} = 1(0)^{2(i_0-1)}1 .
\end{equation*}
 In the admissible words, after 1, the length of the longest sequence of zeros is $ 2(i_0-1)$. 
\begin{equation*}
1(0)^{2i_0-1}\prec 1(0)^{2(i_0-1)}1. 
\end{equation*}
It follows that for all $w\in L_{\beta}$, $ 1w(0)^{2i_0-1} \notin L_{\beta}$ ($1\in L_{\beta}$ and $ (0)^n \in L_{\beta}$ for all $ n \in \N^*$). Thus, $ S_{-\beta}$ is not transitive for $ \beta < \gamma_0$.

\item Assume that $ d(l_\beta, -\beta) = \overline{ d_1 d_2 \cdots d_{2n-1} }$. Consider a word $ x_1 x_2 \cdots x_k \in L_\beta$ such that 
\begin{equation*}
 d_1 d_2 \cdots d_{2n-1} x_1 x_2 \cdots x_k \in L_\beta.
\end{equation*}
 Then, 
\begin{equation*}
d_1 \cdots d_{2n-1} d_1 \cdots d_k \preceq d_1 \cdots d_{2n-1} x_1 \cdots x_k.
\end{equation*}
 there exists an integer $ m $, $ 1 \leq m \leq k $ such that $ x_i = d_i $ for $ 1 \leq i < m $ and $ (-1)^{2n-1+m}(d_m - x_m) \leq 0 $. That is $ (-1)^m(x_m - d_m) \leq 0$. In other words,
\begin{equation*}
x_1 x_2 \cdots x_k \preceq d_1 d_2 \cdots d_k.
\end{equation*}
In fact, $ uv \prec uw $ implies that $ v \prec w $ if $ \vert u \vert $ is even and $ w \prec v $ if $ \vert u \vert $ is odd. 

But, $ x_1 \cdots x_k \in L_\beta $. Thus, $ d_1 \cdots d_k \preceq x_1 \cdots x_k $. 
Hence, it follows that $ x_1 x_2 \cdots x_k = d_1 d_2 \cdots d_k $. Then, for all $ X=d_1 d_2 \cdots d_{p-1} j \in L_\beta$, with $ (-1)^p(d_p - j) < 0$ (we have $ d_1 \cdots d_{\vert X \vert} \prec X $) there does not exist $ Y \in L_\beta$ such that $ d_1 d_2 \cdots d_{2n-1} Y X \in L_\beta $. This implies that $ S_{-\beta}$ is not transitive. 
\end{itemize}
\end{preu}

\begin{rem}\label{rf1}
It is well-known that all coded system is transitive (see \cite{MR858689}). That is, a non transitive system cannot be coded. Thus, the previous proposition allows us to observe that for all $ \beta$ such that $ d(-\frac{\beta}{\beta+1}, -\beta) $ is periodic with odd period or $ \beta $ is less than the golden ratio $ \gamma_0$, $ S_{-\beta} $ is not a coded system.
\end{rem}
Any word of the code of the $\beta$-shift (given in \eqref{BS}) accepts at right any word of the language of the $\beta$-shift. We want to construct a code in the signed parameter case $-\beta$ with this property. To do that, if $ d(l_{\beta}, -\beta) = (d_i)_{i \geq 1}$, we start by storing these sequences in two groups: at first those for which for all integer $ i $, $ d_{2i}< d_1$, and secondly, the sequences for which there exists an integer $ i $ such that $ d_{2i} = d_1$.

 \begin{itemize}
\item If for all $ i \geq 1 $, $ d_{2i}< d_1$ then, we observe that all concatenation of words of the type $ d_1\cdots d_{2n-1} $ ($n \in \N^*)$ is admissible. Moreover we can add at right of such a word any sequence starting by $ d_1$. Therefore, $d_1\cdots d_{n-1}j $ (with $0\leq j < d_1$ and $ (-1)^n(d_n-j)<0$) can be extended at right by any admissible word. 

\item Suppose there exists an integer $ i $ such that $ d_{2i} = d_1 $. So, $ d(l_\beta, -\beta)$ is defined thanks to sequences of positive integers $ (n_i)_{i \geq 1} $ (increasing) and $(p_i)_{i \geq 1}$ such that: 
\begin{equation}
 d(l_\beta, -\beta) = d_1 d_2 \cdots d_{2n_1-1}d_1 \cdots d_{p_1} d_{2n_1 + p_1} \cdots d_{2n_2-1}d_1 \cdots d_{p_2} d_{2n_2 + p_2} \cdots.
 \label{6}
\end{equation}
 
In $ d(l_\beta, -\beta)$, $d_{2n_i-1+k} = d_k $ for all integer $ k $ satisfying $ 1 \leq k \leq p_i $. If $ p_i = 2n_i-1 $, $ d(l_\beta, -\beta)$ is periodic with odd period. If $(d_i)_{i \geq 1} $ is not periodic with odd period, $ p_i $ satisfies both following conditions: $ p_i < 2n_i-1 $ and $ (-1)^{p_i+1}(d_{p_i+1}-d_{2n_i+p_i})<0$ since $ d_1 \cdots d_{p_i}d_{2n_i+p_i}$ is an admissible word.  

Note that if \eqref{6} is satisfied, all concatenation of words $ d_1 \cdots d_{2n-1}$ is no longer admissible like in the previous item. 

We assume $ (d_i)_{i\geq 1}$ non periodic with odd period and we set 
\[B_{i} = d_1 \cdots d_{2n_i-1}. \] 
\begin{rem}\label{rem2}
The word $ X = B_{k_1} \cdots B_{k_t} $ over $\{0,1, \cdots, d_1\} $ is admissible if and only if $ p_{k_i} \leq 2n_{k_{i+1}}-1$, with $ 1 \leq i \leq t-1 $. 
\end{rem}

The words $ d_1 \cdots d_{2k-1} $ such that
\begin{equation}
 \text{ $ 2n_i+p_i \leq 2k-1 \leq 2n_{i+1}-3 $ } \label{Flo1}
 \end{equation}
 for some integer $ i $ (we suppose that $ 2n_0+p_0 = 0 $), and  
\begin{equation}
B_{k_1}B_{k_2} \cdots B_{k_m} d_1\cdots d_{2k-1} \label{Flo2}
\end{equation}
can be extended at right by any sequence starting by $ d_1$ when
 $ p_{k_i} \leq 2n_{k_{i+1}}-1$ for $ 1\leq i \leq m-1$, $ 2k-1 $ satisfying \eqref{Flo1} and $ 2k-1 > p_{k_m} $. 
We set 
\begin{align}
\Delta_{odd}^0 & = \{ d_1 \cdots d_{2n-1}, 2n_i+p_i \leq 2n-1< 2n_{i+1}-1 \vert i \in \N^{*} \text{ or } n < n_1 \} \label{Dood0}, \\
\Delta_{odd}^1 & = \{ B_{k_1}\cdots B_{k_m} X\vert p_{k_i}< 2n_{k_{i+1}}-1, X \in \Delta_{odd}^0, \vert X \vert > p_{k_m} \} \label{Dood1},
\end{align}

Moreover, if we want the word $d_1 \cdots d_{n-1}j $ such that for all $ x \in L_{\beta}$, $d_1 \cdots d_{n-1}jx \in L_{\beta} $, it is necessary to require the following conditions on $ j $:
\begin{equation}
\begin{cases} 
   (-1)^n(d_n-j)< 0, \hspace{0.2cm} 0 \leq j < d_1  \\
  2n_i+p_i + 1 \leq n \leq 2n_{i+1}-1, \hspace{0.2cm} i \in \N 
   \end{cases}
	\label{Flo3}
\end{equation}
with $ 2n_0+p_0 = 0 $. If $ n = 2n_i+p_i $ for some positive integer $ i $, \[ d_1 \cdots d_{n-1}j = d_1 \cdots d_{2n_i-1} d_1 \cdots d_{p_i}j;\] the admissibility of this word implies that:
\begin{equation*}
d_1 \cdots d_{2n_i-1} d_1 \cdots d_{p_i}d_{2n_i+p_i} \preceq d_1 \cdots d_{2n_i-1} d_1 \cdots d_{p_i}j 
\end{equation*}
and
\begin{equation*}
d_1\cdots d_{p_i}d_{p_i+1} \prec d_1 \cdots d_{p_i}j,
\end{equation*}
and thus: 
\begin{equation}
\begin{cases} 
   (-1)^{2n_i+p_i}(d_{2n_i+p_i}-j)< 0, \hspace{0.2cm} 0\leq j < d_1  \\
   (-1)^{p_i+1}(d_{p_i+1}-j)<0.  
   \end{cases}
\end{equation}
That is
\begin{equation}
 (-1)^{p_i}d_{p_i+1} > (-1)^{p_i}j > (-1)^{p_i}d_{2n_i+p_i} \label{eq2nf}.
\end{equation}
So, we define the sets $ \Gamma_0$, $ \Gamma_0^{'}$, $\Gamma_1$ and $ \Gamma_1^{'}$ as follows:
\begin{equation}
x \in \Gamma_0\Leftrightarrow \begin{cases}
                              x=d_1 \cdots d_{n}j, &\text{ with $n \in \N $}\\
															(-1)^{n+1}(d_{n+1}-j)< 0, &\text{ with $ 0\leq j<d_1 $},\\
															 2n_i+p_i\leq n \leq 2n_{i+1}-2, &\text{ $i\in \N$, 
															$n_0=p_0 = 0$}. 
															\end{cases}
\end{equation}
In $ \Gamma_0^{'} $, we have words of the type $ d_1\cdots d_{2n_i+p_i-1} j $ with $ j $ satisfying \eqref{eq2nf}.
\begin{equation}
x\in \Gamma_0^{'} \Leftrightarrow \begin{cases}
                                  x= d_1\cdots d_{2n_i+p_i-1} j, &\text{ with $i \in \N^*$},\\
																	(-1)^{p_i}d_{p_i+1}>(-1)^{p_i}j>(-1)^{p_i}d_{2n_i+p_i};
																	\end{cases}
\end{equation}
\begin{equation}
x \in \Gamma_1 \Leftrightarrow \begin{cases}
                    x=B_{k_1}\cdots B_{k_m}y &\text{ with $y\in\Gamma_0$, $|y|\geq p_{k_m}+2$ };\\                    k_1, \cdots, k_m \in \N^*\\  
		p_{k_i}<2n_{k_{i+1}}-1,  &\text{ $ 1\leq i \leq m-1 $ }.
													     \end{cases}
\end{equation}
We denote by $ \Gamma_1^{'} $ the set of admissible words of the form\[ B_{k_1} \cdots B_{k_{m-1}} B_{k_m}d_1 \cdots d_{p_{k_m}}j, \]
 with $ j $ satisfying \eqref{eq2nf} and for $ 1 \leq i \leq m-1$, $ p_{k_i}< 2n_{k_{i+1}}-1$. In fact, $ B_{k_m}d_1 \cdots d_{p_{k_m}}j $ is just a word of $ \Gamma_0^{'} $ having a length greater than $ p_{k_{m-1}}$. So,
\begin{equation}
x\in \Gamma_1^{'} \Leftrightarrow \begin{cases}
x=  B_{k_1}\cdots B_{k_{m-1}} y\\
y\in \Gamma_0^{'} &\text{ with $ |y|\geq p_{k_{m-1}}+1$},\\
 p_{k_i} < 2n_{k_{i+1}}-1 &\text{ with } 1\leq i \leq m-2.
\end{cases}
\end{equation}
If $(d_i)_{i \geq 1}$ is periodic with odd period, we use the sequence $(d_i^*)_{i\geq 1}$ in the definition of $ \Delta_{odd}^0 $, $ \Delta_{odd}^1 $, $ \Gamma_0 $, $ \Gamma_0^{'} $, $ \Gamma_1$  and $ \Gamma_1^{'} $ instead of $(d_i)_{i \geq 1}$.
\end{itemize}
We set
\begin{equation}
\Delta_{odd} = \begin{cases}
\Delta_{odd}^0\cup \Delta_{odd}^1 &\text{ if \eqref{6} is satisfied }\\
\{d_1\cdots d_{2k+1} | k \in \N \} & \text{ if $ d_{2i}< d_1 $, $ \forall i \in \N^* $}
\end{cases}
\label{Dood},
\end{equation}
and 
\begin{equation}
\Gamma = \begin{cases}
\Gamma_0 \cup \Gamma_0^{'} \cup \Gamma_1 \cup \Gamma_1^{'} &\text{ if \eqref{6} holds} \\
\{d_1 \cdots d_{n-1}j | (-1)^n(d_n-j)<0, 0 \leq j < d_1, n\in \N^* \} & \text{ otherwise}.
\end{cases}
\label{G}
\end{equation}
In \eqref{G}, if $n=1$, $ d_1 \cdots d_{n-1} = \varepsilon $ (the empty word) and thus, $ d_1 \cdots d_{n-1} j =j$. 
\begin{exple}
If $ d_1d_2\cdots = 302\overline{1}$, we have 
\begin{equation*} 
\begin{aligned}
 \Delta_{odd} &= \{3, 302, 30211, 3021111, \cdots \}\\
 \Gamma &= \{0,1,2, 31, 32, 300, 301, 3022, 30210, 302112, \cdots \}.
\end{aligned}
\end{equation*}
\end{exple}
\begin{exple}
Let $ \beta $ be the algebraic integer satisfying $ \beta^4+2 \beta^3+\beta^2-\beta - 1 = 0 $;
\begin{equation*}
d(l_\beta, -\beta) = 2012\overline{1}. 
\end{equation*}
The sequence $ (n_i)_{i \geq 1} $ is finite: $ 2n_1 - 1 = 3 $, $ p_1 = 1 $,  and thus we have $ d_1 \cdots d_{2n_1-1} = 2 0 1 $, $ d_1\cdots d_{p_1} = 2$, and $ d_{2n_1 + p_1} = 1 $. 
\end{exple}
Now, we can give a language of admissible words with properties similar to those of the code of the $\beta$-shift. 

At right of a word of $ \Delta_{odd} $, we can add any admissible word starting by $ d_1$.Thus the free monoid $ \Delta_{odd}^* $ generated by $\Delta_{odd}$ is a subset of $ L_{\beta}$. Moreover all concatenation of a word of $ \Delta_{odd}^*$ and a word of $ \Gamma$ starting by $ d_1 $ is admissible. Let $ \C $ be the language defined by:
\begin{equation}
\C = \{ x y | x \in \Delta_{odd}^*, y \in \Gamma, |y|\geq 2 \} \cup \Gamma  \label{C}.
\end{equation}
By definition of $ \Gamma$, it is obvious to see that at right of each element of $ \C $, we can add  any admissible word. 

If $ \beta $ is less than or equal to the golden ratio $\gamma_0$, $ \C = \{ 0 \} $. However if $ \beta > \gamma_0 $, by construction, $ \C $ is a prefix code on $ \A = \{0, 1, \cdots, d_1 \}$. 

Let $ \Delta_{evn}^0 $ and $ \Delta_{evn}^1 $ be the sets defined as follows:
\begin{align}
\Delta_{evn}^0&= \{ d_1 \cdots d_{2n} \vert n \in \N \} \label{Devn0}, \\
\Delta_{evn}^1&= \{ B_{k_1}\cdots B_{k_m} X\vert p_{k_i}< 2n_{k_{i+1}}-1, X \in \Delta_{evn}^0 \} \label{Devn1},
\end{align}
When $ n = 0 $, $ d_1 \cdots d_{2n}$ is the empty word $ \varepsilon $. We set
\begin{equation} 
\Delta_{evn}  = \begin{cases}
\Delta_{evn}^0\cup \Delta_{evn}^1 &\text{ if \eqref{6} is satisfied } \\
\{ d_1\cdots d_{2n} \vert n \in \N \} & \text{ if $ d_{2i}< d_1 $, $ \forall i \in \N^* $}
\end{cases}
\label{Devn}.
\end{equation}

\begin{lem}\label{corollary 1}
For all $ (k_1, \cdots, k_t) \in \N^{*t}$, $ 1 \leq n \leq p_{k_t}$, $p_{k_i} < 2n_{k_{i+1}}-1$ with $ 1 \leq i \leq t-1$, 
\begin{equation*}
B_{k_1}\cdots B_{k_t}d_1 \cdots d_n \in \Delta_{evn}.
\end{equation*}
\end{lem}

\begin{preu}
It is enough to see that for $ 1 \leq n \leq p_i $, for all $ i \in \N^{*}$, 
\begin{equation*}
d_1\cdots d_{2n_i-1+n} = d_1 \cdots d_{2n_i-1}d_1\cdots d_n.
\end{equation*}
If $ n $ is odd, $ 2n_i-1+n $ is even and then, $ d_1 \cdots d_{2n_i-1+n} \in \Delta_{evn}^{0} \subset \Delta_{evn}$. 

If $ n $ is even, $ d_1 \cdots d_{2n_i-1+n} $ can be seen as a concatenation of $ d_1\cdots d_{2n_i-1} $ and $ d_1 \cdots d_n $. So, 
\begin{equation*}
d_1 \cdots d_{2n_i-1+n} = d_1 \cdots d_{2n_i-1}d_1\cdots d_n \in \Delta_{evn}^1 \subset \Delta_{evn}.
\end{equation*} 
\end{preu}\\
From Lemma \ref{corollary 1}, we can see $ \Delta_{evn} $ as the set of admissible concatenations of words of the type $ B_i $ eventually extended at right by $ d_1 \cdots d_n $ with $ n $ even. So, if $x$ is an admissible word:
\begin{itemize}
\item $ x $ begins by a word of $ \C $, or
\item $ x $ is an admissible concatenation of words of the type $ d_1\cdots d_{2k+1} $ eventually extended at right by $ d_1\cdots d_{2n}$.
\end{itemize} 
If we set 
\begin{equation}
\begin{aligned}
D&=\left\lbrace xy \vert x\in \Delta_{odd}^*\cup \{\varepsilon\}, \hspace{0.1cm} y \in \Delta_{evn} \right\rbrace \label{d} \\
 &=\left\lbrace x y | x \in \Delta_{odd}, \hspace{0.1cm} y \in D \right\rbrace \cup \Delta_{evn},
\end{aligned}
\end{equation}
then, the language $ \la $ is given by:
\begin{equation}
\begin{aligned}
\la = &\left\lbrace u v \vert u \in \C, v \in \la \right\rbrace \cup D.  \label{la}
\end{aligned}
\end{equation}

\begin{rem}\label{r1}
From Theorem 3 of \cite{NguemaNdong20161}, if $ \beta> \gamma_0$ (that is $(d_i)_{i \geq 1} \prec 1\overline{0}$), all sequence between $ d_1\overline{(d_1-1)0} $ and $ \overline{(d_1-1)0} $ cannot be an expansion of $-\frac{\beta}{\beta+1} $ for some $ \beta> 1$ except $\overline{d_1}$. Then, there exists $ y \in \Gamma_1$ such that $ \vert y \vert \geq 2$. Indeed, 
\begin{equation*}
\left[ (d_i)_{i \geq 1} = d(l_\beta, -\beta)\text{ and } \beta> \gamma_0 \right] \Rightarrow (d_i)_{i \geq 1} \prec d_1\overline{(d_1-1)0}.
\end{equation*}
So, we can find $ n\in \N  $ such that $ d_1\overline{(d_1-1)0}^n $ or $ d_1\overline{(d_1-1)0}^n(d_1-1) $ belongs to $\Gamma_1$.
\end{rem}

\begin{theo}\label{p2}
Let $ \beta> \gamma_0 $ and $ d(-\frac{\beta}{\beta+1}, -\beta) = (d_i)_{i \geq 1}$. We assume that $(d_i)_{i \geq 1}$ is not periodic with odd period. Then for all $ n $, $ d_1 \cdots d_n \in L_{\C^*}$.
\end{theo}

\begin{lem}\label{eflem}
Let $ \beta > 1 $. We assume that $(d_i)_{i\geq 1} =  d(l_{\beta}, -\beta) $ is not periodic with odd period and it satisfies \eqref{6}. If there exists an integer $ i_0 $ such that for all $ t \geq i_0 $, $ d_{2n_t} \cdots d_{2n_{t+1}-1} \preceq d_{2n_{i_0}} \cdots d_{2n_{i_0+1}-1} $, then
\begin{equation}
d_1 \cdots d_{2n_{i_0}-1} \overline{ d_{2n_{i_0}} \cdots d_{2n_{i_0+1}-1}} \preceq (d_i)_{i \geq 1} \prec \overline{ d_{2n_{i_0}} \cdots d_{2n_{i_0+1}-1}}.
\end{equation}
\end{lem}

\begin{preu}
We set $ u = d_1 \cdots d_{2n_{i_0}-1} $ and $ v = d_{2n_{i_0}}\cdots d_{2n_{i_0+1}-1} $. Since $ (d_i)_{i \geq 1} $ is not periodic with odd period, we have $ u \prec v $. Since $(d_i)_{i \geq 1} $ starts by $ u $, it becomes obvious that $ (d_i)_{i \geq 1} \prec \overline{v}$.

If $(d_i)_{\geq 1} \neq u \overline{v} $, $ (d_{i})_{i \geq 2n_{i_0}} \neq \overline{v}$. There exists a non negative integer $ m $ such that 
\begin{equation}
(d_i)_{i \geq 2n_{i_0}} = \begin{cases}
 v^m d_{2n_{t}}\cdots d_{2n_{t+1}-1}d_{2n_{t+1}}\cdots d_{2n_{t+2}-1} \cdots  \\
d_{2n_{t}}\cdots d_{2n_{t+1}-1} \prec v. 
\end{cases}
\end{equation}
Since the length of $ v $ is even, it follows that $ (d_i)_{i \geq 2n_{i_0}}\prec \overline{v}$. We obtain the result by adding at left of both words $ u = d_1 \cdots d_{2n_{i_0}-1} $ which is of odd length $ 2n_{i_0}-1 $ and using the property of the alternating order given in \eqref{ef1}. 
\end{preu}

\begin{preu 5}
From Remark \ref{r1}, if $ \beta >   \gamma_0$, $ \Gamma_1 $ contains at least one word $ y $ such that $ \vert y \vert \geq 2 $.
 
If $ d_{2i} < d_1 $, $ d_1 \cdots d_{2n-1} y \in \C$ for all $ n \in \N^{*}$ and $ y \in \Gamma $ with $ \vert y \vert \geq 2$. Thus, for all $ n \in \N^{*} $, $ d_1 \cdots d_n $ is the beginning of a word of $ \C $.

 We assume \eqref{6} satisfied and $(d_i)_{i\geq 1} $ is not periodic with odd period. Suppose  $ d_1\cdots d_n \notin L_{\C^*} $. 

Let $k_0$ be the smallest integer such that $ d_1\cdots d_{k_0} \notin L_{\C^*} $ and $ i> 0$ such that $ 2n_{i}-1 \leq k_0 < 2n_{i+1}-1$. 
Note that for all $ n \geq k_0$, $ d_1 \cdots d_n \notin L_{\C^*}$. Then, for all $ t \geq i $, $ d_1 \cdots d_{p_t}d_{2n_t+p_t}$ and $ d_1 \cdots d_{p_t+1}$ are consecutive (with respect to the alternating order). It follows that 
\begin{equation}
d_{2n_t+p_t} = d_{p_t+1}-(-1)^{p_t} \label{ef2},
\end{equation} 
otherwise, $ d_1 \cdots d_{2n_t+p_t-1}(d_{p_t+1}-(-1)^{p_t}) \in \C$.
Also, there does not exist an integer $ k $ such that $ 2n_t+p_t \leq 2k-1 < 2n_{t+1}-1 $, otherwise $ d_1 \cdots d_{2k-1} y \in \C $ for all $ y \in \C $ with $ \vert y \vert \geq 2 $. Therefore 
\begin{equation*}
2n_{t+1}-1 = \begin{cases}
         2n_t+p_t   &\text{ if $ p_t $ is odd }\\
         2n_t+p_t+1 &\text{ if $ p_t $ is even}.
         \end{cases}
\end{equation*}
If $ p_t$ is even, $ d_{2n_t+p_t+1} = 0 $, otherwise
$ d_1\cdots d_{2n_t-1}d_1\cdots d_{p_t}d_{2n_t+p_t}0 \in \Gamma_1 \subset \C$.
So
\begin{equation}
d_1\cdots d_{2n_{t+1}-1}= \begin{cases} 
                        d_1\cdots d_{2n_t-1}d_1\cdots d_{p_t}(d_{p_t+1}+1) &\text{ if $ p_t $ is odd }\\
                        d_1\cdots d_{2n_t-1}d_1\cdots d_{p_t}(d_{p_t+1}-1)0 &\text{ if $ p_t $ is even }.
                          \end{cases} \label{ef3}
\end{equation}
Observe that if $ d_1 \cdots d_{p_t + 1}$ is the beginning of a word of $ \C $, automatically $ d_1 \cdots d_{2n_t-1} d_1 \cdots d_{p_t+1}$ is the beginning of a word of $\C $. Thus, $ d_1 \cdots d_{2n_t-1} \not\in L_{\C^*}$ implies that $ d_1 \cdots d_{p_t+1} \not\in L_{\C^*}$. Since $(d_i)_{i \geq 1}$ is not periodic with odd period, $ p_i < 2n_i-1$. But, $ k_0$ is supposed to be greater than $ 2n_i - 1$. Thus we have $p_i < k_0$, and then for any integer $ t \geq i $, $ p_i < k_0 \leq p_t +1$. That is $ p_i \leq p_t $. In particular $ p_i < 2n_i-1 \leq k_0 \leq p_i+1$ (and then, $ p_i+1 = 2n_{i}-1$).

If $ p_t = p_i $, \eqref{ef3} requires $ d_{2n_i} \cdots d_{2n_{i+1}-1}= d_{2n_t} \cdots d_{2n_{t+1}-1} $. 
If $ p_t> p_i $, we have 
\[ d_{2n_t} \cdots d_{2n_{t+1}-1} \prec d_{2n_i} \cdots d_{2n_{i+1}-1}. \] 
Then,
\begin{equation*}
d_1\cdots d_{p_i}d_{p_i+1}d_{p_i+2} = 
                           d_1\cdots d_{2n_i-2}d_{2n_i-1} d_1. 
\end{equation*}
We set 
\begin{align*}
U_0&= d_1\cdots d_{2n_i-1} \\
 V_0& = d_{2n_i} \cdots d_{2n_{i+1}-1} = d_1 \cdots d_{p_i}d_{2n_i+p_i} \cdots d_{2n_{i+1}-1}.
\end{align*}
 $ V_0 $ is such that:
\begin{equation}\label{eq16}
V_0 =
                 d_1 \cdots d_{2n_i-2}(d_{2n_i-1}-1)0. 
\end{equation}
From Lemma \ref{eflem}, $ U_0 \overline{V_0}\preceq (d_i)_{i \geq 1} \prec \overline{V_0}$.
 From Proposition 9 and Theorem 3 of \cite{NguemaNdong20161}, the unique sequence between $ U_0\overline{V_0} $ and $ \overline{V_0}$ which is the $(-\beta)$-expansion of $l_{\beta}$, for some $ \beta > \gamma_0$ is $ \overline{U_0}$. This is absurd, since $ (d_i)_{i \geq 1 } = d(l_{\beta}, -\beta)$ is supposed to be non-periodic with odd period and $ (d_i)_{i \geq 1} \not\neq \overline{V_0}$ because $(d_i)_{i \geq 1} $ starts by $ U_0$. Then the assumption $ n \geq 2n_i-1$, $ d_1 \cdots d_n $ is not the beginning of a word of $ \C $ is inaccurate. 
\end{preu 5}

From Proposition 9 of \cite{NguemaNdong20161}, when $ \overline{V_0} \prec 1\overline{0} $, the sub-shift of infinite words for which all sub-words is bigger than $ U_0 \overline{V_0} $ and the sub-shift of infinite words for which all sub-words is bigger than $ \overline{V_0} $ have the same entropy. We denote by $ \psi $ the map from $ \{U_0, V_0 \} $ into $ \{U_0, V_0 \}^* $ defined by $ \psi(U_0) = U_0V_0 = U_1 $ and $ \psi(V_0) = U_0U_0 = V_1 $. 
The limit of non-periodic infinite words $(d_i)_{i \geq 1} $ such that $ U_0 \overline{V_0} \preceq (d_i)_{i \geq 1} \preceq \overline{V_0} $ and $ (d_i)_{i\geq 1} \prec (d_{i})_{i \geq k } $ for all $ k > 1 $ (not necessary the $(-\beta)$-expansion of $ l_{\beta} $ for some $ \beta> \gamma_0$) is $ \psi^{\infty}(U_0) = \lim\limits_{n \mapsto \infty } \psi^n(U_0) $.
\begin{equation}
\psi^{\infty}(U_0) = U_0V_0U_0U_0U_0V_0U_0V_0U_0V_0U_0U_0U_0V_0U_0U_0\cdots.
\end{equation}
In fact, if $ \psi^{\infty}(U_0) \prec (d_i)_{i \geq 1} \prec \overline{V_0}$, then there exists $ n \in\N^* $ such that $ (d_i)_{i \geq 1} = \overline{\psi^n(U_0)} $.

If $ U_0 = d_1 \cdots d_{2n_i-1} = 1 $, $ V_0 = 00 $, $ \beta \leq \gamma_0$. 
\begin{equation*}
\psi^{\infty}(1) = 1001110010010011100111001110010010011100100\cdots.
\end{equation*}
In this case, the $(-\beta)$-expansion of $ l_{\beta}$ is between $U_0 \overline{V_0} = 1 \overline{00}$ and $ \psi^{\infty}(1)$. 
\begin{equation*}
\psi^{\infty}(1) = \lim\limits_{ \beta \mapsto 1} d(l_{\beta}, -\beta). 
\end{equation*}

\begin{lem}\label{l2}
Let $ \beta $ be a real number bigger than the golden ratio $ \gamma_0$ and $ d(l_\beta, -\beta) $ non-periodic with odd period. Then, the $(-\beta)$-shift $ S_{-\beta} $ is coded by $ \C $.  
\end{lem}
\begin{preu}
Suppose $ \beta > \gamma_0$. Then, $ \C $ is a prefix code. From Theorem \ref{p2}, for all $ n \in \ N^* $, $ d_1 \cdots d_n \in L_{\C^*}$. It follows that $ D \subset L_{\C^*} $. To conclude, it suffices to observe that (since the empty word $\varepsilon $ belongs to $ D $):
\begin{equation*}
\{ x y | x \in \C, y \in \la \} =  \{ x y | x \in \C^{*} \text{ and $ y \in D$} \}.
\end{equation*}
\end{preu}

From the proof of Theorem \ref{p2}, if $ \beta > \gamma_0 $ and $ d(l_{\beta}, -\beta) = \overline{d_1 \cdots d_{2n_i-1} }$, we have $ d_{2n_i-1}\neq 0 $ and the word $ d_1\cdots d_{2n_i-1}$ is an intransitive. The $(-\beta)$-shift, as defined in \eqref{S}, is not coded, but contains a coded sub-shift: the dynamical system of words for which, in the meaning of alternating order, all sub-word is bigger than $ \overline{d_1\cdots d_{2n_i-2}(d_{2n_i-1}-1)0}$. Moreover, both systems have the same entropy (see \cite{NguemaNdong20161}).

If we consider the definition of the $ (-\beta)$-shift given in the introduction in the case where $ \beta $ is integer, that is, the expansion of $ -\frac{\beta}{\beta+1} $ is periodic with period 1, we obtain  
\begin{equation*}
 S_{-\beta} = \left\lbrace (x_i)_{i \in \Z}; \overline{\beta} \preceq (x_i)_{i \geq n} \preceq \overline{(0, \beta-1)}, \forall n \right\rbrace. 
\end{equation*}
But in this special case, the definition given by S. Ito and T. Sadahiro is
\begin{equation}
 S_{-\beta} = \left\lbrace (x_i)_{i \in \Z};  \overline{(\beta-1,0)} \preceq (x_i)_{i \geq n} , \forall n \right\rbrace.
 \label{(24)}
\end{equation}
 In fact, all sequence $ (x_i)_{i \in \Z} $ of $ S_{-\beta} $ satisfies, for all $ n $
\begin{equation}
  \lim\limits_{ y \rightarrow l_{\beta}^+}d(y, -\beta) \preceq x_nx_{n+1}\cdots \preceq  \lim\limits_{ y \rightarrow r_{\beta}^-}d(y, -\beta) 
  \label{(25)}
\end{equation}
with $ l_{\beta} = -\frac{\beta}{\beta+1} $ and $ r_{\beta} = \frac{1}{\beta+1} $. Using Lemma 6 of \cite{MR2534912}, $ \lim\limits_{ y \rightarrow l_{\beta}^+}d(y, -\beta)=  \overline{(\beta-1,0)} $ and $ \lim\limits_{ y \rightarrow r_{\beta}^-}d(y, -\beta) = \overline{(0, \beta-1)} $. In this form, $ S_{-\beta} $ is coded. Generally, the symbolic dynamical system $ S_{-\beta} $ contains a sub-shift coded by a prefix code and which is the support of the maximal entropy measure. When $ d(l_\beta, -\beta)$ is periodic with odd period $ 2n-1$, we introduce the sub-shift $\tilde{S}_{-\beta} $ (corrected $(-\beta)$-shift) defined by:
 \begin{defi}
\begin{equation*}
 \tilde{S}_{-\beta} = \left\lbrace (x_i)_{i \in \Z} ; d_1^*d_2^*\cdots \preceq x_k x_{k+1} \cdots, \text{ $ \forall k $ } \right\rbrace 
\end{equation*}
with 
 \begin{equation*}
 (d_{i}^{*})_{i \geq 1} = \begin{cases}
                         \overline{(d_1, \cdots, d_{2n-2}, d_{2n-1}-1, 0)} &\text{ if $ (d_i)_{i \geq 1} = \overline{(d_1, \cdots, d_{2n-1} )} $ }\\
                         (d_i)_{i \geq 1} &\text{ otherwise }. 
                         \end{cases}                        
\end{equation*} 
\end{defi} 
Using Lemma 6 of \cite{MR2534912}, we see easily that 
$ (d_1^*,d_2^*, \cdots )= \lim\limits_{ x \rightarrow l_{\beta}^+}d(x, -\beta)$ (corrected $(-\beta)$-expansion of $ l_{\beta}$).

All real has a representation in $\tilde{S}_{-\beta} $ since $\sum\limits_{i\geq 1}\frac{d_i^{*}}{(-\beta)^i} = l_\beta$ (see the proof of Proposition 8 of \cite{NguemaNdong20161}). 
We find more convenient to use $\s$ as $(-\beta)$-shift instead of $ S_{-\beta}$.
In fact, the sequence $(d_i^*)_{i \geq 1}$ plays the role of $0.9999999\cdots $ in base 10. For instance, $0.9999999\cdots $ is the representation of 1 in base 10. 
\begin{equation*}
1= 0.9999999999\cdots
\end{equation*}

Now, we are ready to yield the proof of Theorem \ref{t1}.

\begin{preu 2}
When $ \beta = \gamma_0$, it is easy to see that the system is coded by 
\[\Delta_{odd} = \{1, 100, 10000, \cdots \}.\]
 But this code is not optimal. Indeed, this language is obtained owing to two words: 1 and 00. Thus, we code $S_{\gamma_0}$ by $ \{1, 00\}$. To complete the proof of Theorem \ref{t1}, it is enough to use Lemma \ref{l2} and Proposition \ref{prop1}.
\end{preu 2}

 From Theorem \ref{t1}, $\tilde{S}_{-\beta}$ is a coded (by $\C $ defined in \eqref{C}) if only if $ \beta \geq \gamma_0$.
Throughout the rest of this paper, we focus our interest in the study of $ \tilde{S}_{-\beta}$ instead of $ S_{-\beta}$. 

\subsection{Recurrent positive code}  

 We have seen in the previous subsection that when $ \beta > \gamma_0 $, the symbolic dynamical system $ \s $ is coded by the language $ \C $ and $ D $ is a subset of $ L_{\C^*} $. In fact, $ \s $ can be seen as the support of the maximal entropy measure. If furthermore $ d_{2i} < d_1 $ for all integer $ i > 0$, $ \C $ and $ \Delta_{odd} $ allow us to characterize words of $ \la $. 

When \eqref{6} is satisfied, it becomes increasingly unclear because of the non admissibility of certain concatenations of $ B_i $ (see Remark \ref{rem2}). 
In particular, for $ \beta < \gamma_0$, we know that $\s$ (and then $ S_{-\beta}$) is not transitive and $ \C = \{0\}$. This implies that the support of the maximal entropy measure is included in $ D$. It may asked: what is this support? 
With a view to determining this support, we study in detail hereafter the conditions of admissibility of product of $ B_i $ and also, we get another formulation of the set $D$. 
We exhibit different codes for the writing of these concatenations. If a language $ \Omega $ is one of these codes, any concatenation in $ \Omega $ is admissible. 

\begin{rem}\label{rem4}
Let $ X=B_{k_1}\cdots B_{k_m} $ be an admissible word of a code, $ X X \in L_{\beta}$. From Remark \ref{rem2}, $ p_{k_i} \leq 2n_{k_{i+1}}-1$ with $ 1\leq i \leq m-1$, and $ p_{k_m} \leq 2n_{k_1}-1$. If $ 2n_i-1 < p_{t_m} \leq 2n_{i+1}$ for some $ i $, $ X $ can be extended at right by any word starting by $ B_{i+1}$. 
\end{rem}
We set 
\begin{equation}
J(0) = \{ t, p_t \leq 2n_1-1 \}, \label{J0}
\end{equation}
and for all $ i $, 
\begin{equation}
J(i) = \{ t, 2n_i-1 \leq p_t < 2n_{i+1}-1 \} \label{Ji}.
\end{equation}
Let $ \Delta^{(i)}$ be the sets such that:
\begin{equation}
x \in \Delta^{(i)} \Longleftrightarrow 
                    \begin{cases} 
                    x = B_{t_1}\cdots B_{t_m},                      \\
                    p_{t_k} \leq 2n_{t_{k+1}}-1, p_{t_m}< 2n_{t_1}-1 \\
                     t_m \in J(i),  \\
                      t_k \in J(l) \text{ for $ k\neq m $ and $ l \geq i+1$ }.
                    \end{cases} \label{Di}
\end{equation} 
Let us explain a little bit the definition of the set $ \Delta^{(i)}$. Consider an element $ x =B_{t_1}\cdots B_{t_m}$ of $ \Di $. 
\begin{itemize}
\item The condition $ t_k \in J(l) $ with $ l \geq i+1$ involves that $ t_k \not\in J(i)$ and then $ x $ cannot be a concatenation of words of $\Di$. That is $ \Di $ is a (prefix or suffix) code.
\item $ p_{t_k} \leq 2n_{t_{k+1}}-1 $ allows to have $ x $ admissible, but when $p_{t_m}< 2n_{t_1}-1$, the word generates a periodic expansion.
\item For the condition $ l \geq i + 1$, if we suppose $ t_1 \in J(l) $ with $ l < i $, the word $ B_{t_2}\cdots B_{t_m}B_{t_1} \in \Delta^{(l)}$ and it is a result of a permutation of $ x $. In fact, if $ t_1 \in J(l) $ with $ l < i $, $ x $ is a word of the language of the free monoid generated by $ \Delta^{(l)}$. Thus $ l \geq i+1$ ensures the fact that two words $ x \in \Di $ and $ y \in \Delta^{(i)}$ cannot generate the same periodic orbit. 
\end{itemize}

From Remark \ref{rem2}, $ (\Delta^{(i)})^* \subset \la $.

\begin{lem}\label{rem5}
Let $ x \in \Delta^{(i)} $ and $ y \in \Delta^{(j)}$, with $ i \leq j $. Then, $ xy \in \la $. 
\end{lem}
\begin{preu}
Let $ x = B_{t_1}\cdots B_{t_s} \in \Delta^{(i)} $ and $ y = B_{t_{s+1}}\cdots B_{s+m}\in \Delta^{(j)} $. From \eqref{Di}, $ 2n_i-1<p_{t_s} < 2n_{i+1}-1$. Since $ n_{t_{s+1}} \geq n_{n_{j+1}} $, we have 
\begin{equation*}
i \leq j \Rightarrow 2n_{t_{s+1}}-1 \geq 2n_{j+1}-1 \geq 2n_{i+1}-1 > p_{t_s}
\end{equation*}
and we find hence that $ xy \in L_{\beta}$.
\end{preu}

Now, we know the different sub-languages necessary to characterize words of $ \la $: $ \C $, $ \Delta_{odd} $, $ \Delta^{(i)}$ for all $ i \in \N^*$ if the corrected $(-\beta)$-expansion of $l_{\beta}$ satisfies \eqref{6}; or $ \C $ and $ \Delta_{odd} $ if \eqref{6} does not hold.

Let $ \beta $ be a real number bigger than 1. In the previous paragraph, we constructed a prefix code of the $(-\beta)$-shift. But, if a code is recurrent positive, it is more interesting above all if the system is intrinsically ergodic with entropy $ \log \beta $.
Start by giving the following definition:
\begin{defi}
Let $ X $ be a symbolic dynamical system and $ t = t_1 \cdots t_n $ a word of $L_X$. We denotes by $ [t] $, the set of infinite words $(x_i)_{i \geq 1} $ starting by $ t $. That is:
\begin{equation*}
x_1x_2\cdots x_n= t_1t_2 \cdots t_n.
\end{equation*}
\end{defi}

Consider a symbolic dynamical system $ X $. In fact, the existence of a recurrent positive prefix code $ \C$ implies that there exists a probability measure $ \nu$ on the space $W(\C)$, closure $ \C^{\Z} $ defined as: 
\begin{equation*}
\nu([x]) = \frac{1}{\beta^{\vert x \vert}}, \text{ for all $ x \in \C $},
\end{equation*} 
with $ \log \beta $ the entropy of $ X $ endowed with the shift. We have
$W(\C) = \bigcup\limits_{x \in \C}[x] $; for $(x,y)\in \C^2 $, $ x\neq y$, $ [x] \cap [y] = \varnothing $. Thus, 
\begin{align*}
 \nu (W(\C)) &= \sum\limits_{x \in \C} \nu([x]) \\
             &= \sum\limits_{x \in \C} \frac{1}{\beta^{\vert x \vert}} \\
             &= 1.
\end{align*} 
The entropy $ h_{\nu}$ of the probability $\nu$ is:
\begin{align*}
h_{\nu} &= - \sum\limits_{x\in \C} \nu([x])\log \nu([x]) \\
        &= \left( \sum\limits_{x\in \C} \frac{ \vert x \vert }{\beta^{\vert x \vert }}\right) \log \beta.
\end{align*}  
This expression exists since $ \sum\limits_{x\in \C} \frac{ \vert x \vert }{\beta^{\vert x \vert }} < + \infty $ when the code is recurrent positive. Thus, the maximal entropy measure $ \mu $ is given by:
\begin{equation*}
\mu = \left(\sum\limits_{x\in \C} \frac{ \vert x \vert }{\beta^{\vert x \vert }}\right)^{-1} \nu.
\end{equation*}
See for example \cite{MR858689} and \cite{MR939059} for more explanations.   

 Using the formal power series, the following result establishes a link between $ \la $, $\Delta_{odd}$ and $ \Di $, $ i \geq 1$.
\begin{theo}\label{th2}
Let $ \beta $ be a real number ($\beta> 1$), $ \la $ the language of the corrected $(-\beta)$-shift $\s$, $\C $, $ \Delta_{odd}$ and $ \Delta^{(i)}$ be the sets defined in  \eqref{C}, \eqref{Dood},  and \eqref{Di} respectively. In the meaning of the formal power  series, we have the following relations:
\begin{equation}
\sum\limits_{n \geq 0}z^n= 
( 1+z)( 1-\sum\limits_{x\in \C}z^{\vert x \vert})( 1-\sum\limits_{x \in \Delta_{odd}}z^{\vert x \vert})\prod\limits_{i \geq 1}(1-\sum\limits_{x \in \Delta^{(i)}}z^{\vert x \vert})\sum\limits_{x \in L_{\s}}z^{\vert x \vert}\label{eq30}
\end{equation}
if \eqref{6} occurs, or
\begin{equation*}
\sum\limits_{n \geq 0}z^n= \left( 1+z\right) \left( 1-\sum\limits_{x\in \C}z^{\vert x \vert}\right)                                                           \left( 1-\sum\limits_{x \in \Delta_{odd}}z^{\vert x \vert}\right)\sum\limits_{x \in \la}z^{\vert x \vert} 
\end{equation*}
if $ d_{2i}^* < d_1^* $, $ \forall i \in \N^* $ and 
where $(d_i^*)_{i \geq 1}$ is the corrected $(-\beta)$-expansion of $l_{\beta} = -\frac{\beta}{\beta+1}$.
\end{theo}

\begin{preu}
In the following, $ H_n$, $c_n$, $ a_n $ and $ \delta_n $ denote respectively the number of words of length $ n $ in $ \la $, $ \C $, $ \Delta_{odd}$ and $ D $ (given in \eqref{d}).

We have seen that a word of $ \C $ can be extended at right by any word of $ \la$. Therefore, at right of a word of $ \Delta_{odd} $ can be added any word of $ D $. So, from \eqref{la}
\begin{equation}
\sum\limits_{n \geq 0} H_nz^n = \left( \sum\limits_{n \geq 1}c_nz^n \right) \sum\limits_{n \geq 0} H_nz^n +\sum\limits_{n \geq 0} \delta_n z^n. \label{eq31}
\end{equation}
Let us explain a little bit the equation above. The coefficients of the formal power series $ \left( \sum\limits_{n \geq 1}c_nz^n \right) \sum\limits_{n \geq 0} H_nz^n$ count admissible finite sequence starting by a word of $ \C $. 

We denote by $ \sum\limits_{n \geq 1} b_{i,n}z^n $ the formal power series for which the coefficients count the words of $\Delta_{evn}$ which start by a word of $ \Delta^{(j)}$ with $j\geq i$. So, the coefficients of $ \sum\limits_{n \geq 1} b_{0,n}z^n$ count all words of $ \Delta_{evn}$. 

\begin{equation}
\sum\limits_{n \geq 0} \delta_nz^n = \left(\sum\limits_{n \geq 1} a_n z^n\right)\left( \sum\limits_{n \geq 0} \delta_nz^n \right)+\sum\limits_{n \geq 1} b_{0,n}z^n+1. \label{eq32}
\end{equation}
This equation means that in $ D $, we have admissible strings starting by words of $ \Delta_{odd} $ and those belonging to $\Delta_{evn}$. We set 
\begin{equation*}
\sum\limits_{x\in \Delta^{(i)}}z^{\vert x \vert} = P_i.
\end{equation*}
From Lemma \ref{corollary 1}, Remark \ref{rem4} and Lemma \ref{rem5}
\begin{equation}
\sum\limits_{n \geq 1}b_{0,n}z^n= \begin{cases}
                 \sum\limits_{n \geq 1}z^{2n} +\sum\limits_{n \geq 1}b_{1,n}z^n &\text{ if \eqref{6} occurs }; \\
                 \sum\limits_{n \geq 1} z^{2n} &\text{ if $ d_{2i}^* < d_1^* $, $ \forall  i $ }
                 \end{cases}\label{eq33}
\end{equation}
and for all $i\geq 1$,
\begin{equation}
\sum\limits_{n \geq 1}b_{i,n}z^n= P_{i}\left( \sum\limits_{n \geq 0}z^{2n}+\sum\limits_{n \geq 1}b_{i,n}z^n\right) +\sum\limits_{n \geq 1}b_{i+1,n}z^n\label{eq34}.
\end{equation}
From \eqref{eq32}, \eqref{eq33} and \eqref{eq34}, we have
\begin{equation}
\sum\limits_{n \geq 0}z^{2n}= \begin{cases} 
(1-\sum\limits_{n \geq 1}a_nz^n) \prod\limits_{i \geq 1}(1-P_i)\sum\limits_{n\geq 0}\delta_nz^n &\text{ if \eqref{6} occurs }\\
                              ( 1-\sum\limits_{n \geq 1}a_nz^n)\sum\limits_{n\geq 0}\delta_nz^n &\text{ if $ d_{2i}^* < d_1^* $, $ \forall i \in \N^*$ }
                              \end{cases}  \label{eq35}                     
\end{equation}
Thus, to obtain \eqref{eq30}, it is enough to multiply \eqref{eq35} by $ 1+z $ and use \eqref{eq31}. 
\end{preu}

The following corollaries are the consequences of Theorem \ref{th2} and Lemma 1 of \cite{NguemaNdong20161}. 
\begin{coro}
In the disk $ b(0, \frac{1}{\beta})$ of center 0 and radius $ \frac{1}{\beta} $, we have:
\begin{equation}
1-\sum\limits_{n \geq 1} (d_{n-1}^*-d_n^*)(-z)^n=(1+z)( 1-\sum\limits_{x \in \C }z^{\vert x \vert}) ( 1-\sum\limits_{x \in \Delta_{odd}} z^{\vert x \vert}) \prod\limits_{i \geq 1} ( 1-\sum\limits_{x \in \Delta^{(i)}} z^{\vert x \vert} ) \label{eq39}
\end{equation}
with $\sum\limits_{x \in \Delta^{(i)}} z^{\vert x \vert} = 0$ when $ d_{2i}^*<d_1^*$ for all $ i $.
\end{coro}
\begin{preu}
Let $ H_n $ be the number of words of length $ n $ in $ L_{\beta}$. From \cite{NguemaNdong20161}, we know that the formula for the factor complexity of the corrected $ (-\beta)$-shift (or the language $\la$) is given by:
\begin{equation*}
H_n = \sum\limits_{k=1}^{n}(-1)^k(d_{k-1}^*-d_k^*)H_{n-k} + 1
\end{equation*}
with $ (d_i)_{i \geq 1}=d(l_{\beta}, -\beta)$, $ d_0 = 0 $ and $ H_0 = 1$.
By simple calculus in the sense of power formal series, 
\begin{equation*}
\sum\limits_{x\in L_{\s}}z^{\vert x \vert} = \dfrac{\sum\limits_{n \geq 0}z^n}{1-\sum\limits_{n \geq 1}(-1)^n(d_{n-1}^*-d_n^*)z^n}. 
\end{equation*}
We conclude by using Theorem \ref{th2}.
\end{preu}

\begin{coro}\label{cflo}
The expansion $d(l_{\beta},-\beta)$ is supposed to be periodic with odd period $ 2p-1$. Then, in the sense of formal power series,
\begin{equation}
\sum\limits_{x \in \la}z^{\vert x \vert } = \dfrac{(1-z^{2p})\sum\limits_{n \geq 0}z^n}{1-\sum\limits_{n \geq 1}(-1)^n(d_{n-1}-d_n)z^n} \label{cfl}.
\end{equation}
\end{coro}

\begin{preu}
Any word $(x_i)_{i \geq 1}$ in the $(-\beta)$-shift $ S_{-\beta}$ satisfies
\begin{equation*}
(d_i)_{i \geq 1} \preceq (x_{i + n})_{i\geq 1} \preceq (d_{i-1}^*)_{i \geq 1}, \text{ $ \forall n \in \N $ }.
\end{equation*} 
According to Lemma 1 of \cite{NguemaNdong20161}, the formula for the factor complexity of the $(-\beta)$-shift is: 
\begin{equation}
\tilde{H}_n = \sum\limits_{k=1}^{n}(-1)^k(d_{k-1}^*-d_k)\tilde{H}_{n-k} + 1
\end{equation}
where $ \tilde{H}_n $ denotes the number of words of length $ n $ in the language of the $(-\beta)$-shift. We obtain \eqref{cfl} by using the definition of $(d_i^*)_{i \geq 1} $ given in \eqref{D}.
\end{preu}

From \eqref{eq39} and \eqref{cfl}, $ \sum\limits_{n \geq 0} H_n z^n $ and $ \sum\limits_{n \geq 0} \tilde{H}_n z^n $ have at pole at $ \frac{1}{\beta}$ and it is the smallest pole in modulus. In fact, $ \frac{1}{\beta}$ si the smallest zero in modulus of 
\[ 1-\sum\limits_{k\geq 1}(-1)^k(d_{k-1}^*-d_k)z^k \] 


However, if $ \beta $ is less than the golden ratio, both inclusions can not hold. Indeed, $\C = \{0\}$ and then $ L_{\C^{*}} = \{ \overline{0}^{n}, n \in \N^{*} \} $. The support of the maximal entropy measure is coded by $ \Delta_{odd}$ (in this case $\sum\limits_{x \in \Delta_{odd}}\frac{1}{\beta^{\vert x \vert }}=1$) or by $ \Delta^{(i)}$ for some $ i $ (and $\sum\limits_{x \in \Delta^{(i)}}\frac{1}{\beta^{\vert x \vert }}=1$).

 Recall that the morphism $ \psi$ on $ \{0, 1 \} $ is given by $ \psi (0) =  1 $, $ \psi(1) = 100$ and define the sequences $ (u_n)_{n \geq 0} $ and $ (v_n)_{n \geq 0}$ by $ u_0 = 1 $, $ v_0 = 00$ and for $ n \geq 1$, $u_n= u_{n-1}v_{n-1} = \psi^n(1) $, $ v_n = u_{n-1}u_{n-1}$. From Lemma 2 of \cite{NguemaNdong20161}, $ \vert u_n \vert $ is odd and $ \vert v_n \vert $ is even. In fact, 
 \begin{equation*}
\vert u_n \vert  =  \vert v_n \vert + (-1)^n.
 \end{equation*}
 Moreover, note that there is no word between $ u_k=u_{k-1}u_{k-2}u_{k-2} $ and $ v_k = u_{k-1}u_{k-1}$. Indeed, from the definition of $ u_k $, we observe easily that
\begin{equation*}
u_k = \begin{cases}
      u_{k-1}u_{k-2} u_{k-3}\cdots u_1 u_0 1  &\text{ if $ k $ is even,}\\
      u_{k-1}u_{k-2} u_{k-3}\cdots u_1 u_0 00 &\text{ if $ k $ is odd  }
      \end{cases}
\end{equation*}
and
\begin{equation*}
v_k = \begin{cases}
      u_{k-1}u_{k-2} u_{k-3}\cdots u_1 u_000 &\text{ if $ n $ is even,} \\
      u_{k-1}u_{k-2} u_{k-3}\cdots u_1 u_0 1 &\text{ if $ n $ is odd  }.
      \end{cases}
\end{equation*} 
 
Let $\gamma_n $ be the real number such that
\begin{equation*}
d(l_{\gamma_n}, -\gamma_n) = u_n \overline{v_n}. 
\end{equation*}
$ \gamma_n $ is the largest number satisfying: 
\begin{equation*}
1 = \frac{1}{\gamma_n^{\vert u_n \vert}}+\frac{1}{\gamma_n^{\vert v_n \vert}},
\end{equation*}
that is, $ \gamma_n $ is the largest root of $ X^{l_n} - X - 1 $, where $ l_n = max(\vert u_n \vert, \vert v_n \vert )$. 
The sequence $(\gamma_n)_{n \geq 0} $ decreases until 1 (see Proposition 5 and 6 of \cite{NguemaNdong20161}) and we have $ \lim\limits_{n \rightarrow +\infty}d(l_{\gamma_n}, -\gamma_n) = \psi^{\infty}(1)$. 

The following proposition gives us a writing of $ d(l_{\beta}, -\beta)$ for $ \beta $ less than the golden ratio.
\begin{prop}\label{prop7}
Let $ \beta $ be a real number such that $ 1 < \beta < \gamma_0$. Then, there exists a sequence of integers $ (k_i)_{i \geq 1} $ and $ n \in \N $ such that 
\begin{equation}
d(l_{\beta}, -\beta) = u_n \overline{v_n}^{k_1+1}u_n\overline{v_n}^{k_2}u_n \overline{v_n}^{k_3}\cdots \label{eq36},
\end{equation}
or
\begin{equation}
d(l_{\beta}, -\beta) = u_n v_n (u_nu_n)^{k_1}v_n (u_n)^{k_2}v_n(u_n)^{k_3}v_n \cdots. \label{eq36bis}
\end{equation}
\end{prop}

\begin{preu}
Let $ \beta $ be a real number and suppose $ 1 < \beta < \gamma_0$. Since $ (\gamma_n)_{n \geq 0} $ decreases until 1, there exists an integer $ n $ such that $ \gamma_{n+1} \leq \beta < \gamma_n $. So,
\begin{equation*}
d(l_{\gamma_n}, -\gamma_n) \prec d(l_{\beta}, -\beta) \preceq d(l_{\gamma_{n+1}}, -\gamma_{n+1}).
\end{equation*}
This means
\begin{equation*}
u_n\overline{v_n} \prec d(l_{\beta}, -\beta) \preceq u_{n}v_n\overline{u_n}.
\end{equation*}
It is easy to see that $ u_k \prec v_k $ and there is no word between $ u_k $ and $ v_k $. In an infinite word, $u_k $ is followed by $ v_k$ or by $ u_k $. We obtain \eqref{eq36} by  using the fact that $ u_n\overline{v_n} \prec d(l_{\beta}, -\beta) $. To obtain \eqref{eq36bis}, we interpret the fact that $ d(l_{\beta}, -\beta) \preceq u_{n}v_n\overline{u_n} $.
\end{preu}

S. Ito and T. Sadahiro determined the unique $T_{-\beta}$-invariant measure with maximal entropy on $I_{\beta} = [-\frac{\beta}{\beta+1}, \frac{1}{\beta+1} )$. Note that the structure of the one-side (right) $(-\beta)$-shift endowed with the shift $\sigma$ is transported to $I_{\beta}$ endowed with the $(-\beta)$-transformation. For $ \beta $ taken in the open interval delimited by 1 and the golden ratio, we know that the system is not coded, and then the support of the measure of maximal entropy is a coded subsystem strictly included in the $(-\beta)$-shift. The images by $ T: (x_i)_{i \geq 1} \in S_{-\beta}^r\mapsto \sum\limits_{n\geq 1}\frac{x_i}{(-\beta)^n} \in \overline{I_{\beta}}$ of subsystems non included in the support of the intrinsic ergodic measure correspond to gaps on $ I_{\beta}$. This phenomenon has been closely studied by L. Liao and W. Steiner in \cite{MR2974214}.

 The gaps on $I_{\beta}$ are the intervals $A_{k, i}$ defined as followed:
\begin{equation}\label{eq37}
A_{k, i} = \begin{cases}
           [s_{\vert u_k \vert +i}, s_{\vert u_k \vert + \vert u_{k-1} \vert + i} )  &\text{ if $ i $ is even} \\
           [ s_{\vert u_k \vert + \vert u_{k-1} \vert + i}, s_{\vert u_k \vert +i} ) &\text{ if $ i $ is odd }
           \end{cases}
\end{equation}
with $ k < n $, $ i < \vert u_{k-1} \vert $ and $ s_t = T_{-\beta}^t(l_{\beta})$. Note that 
\begin{equation*}
 u_n = u_ku_{k-1}u_{k-1}u_ku_k \cdots u_{n-2}u_{n-2}. 
\end{equation*}
From \eqref{eq36}, 
\begin{equation*}
d(s_{\vert u_k \vert +i}, -\beta) = \sigma^i(u_{k-1})u_{k-1}u_ku_k \cdots u_{n-2}u_{n-2} \overline{v_n}^{k_1}u_n\overline{v_n}^{k_2}u_n \overline{v_n}^{k_3}\cdots
\end{equation*}
and
\begin{equation*}
d(s_{\vert u_k \vert + \vert u_{k-1}\vert + i}, -\beta) = \sigma^i(u_{k-1})u_ku_k \cdots u_{n-2}u_{n-2} \overline{v_n}^{k_1}u_n\overline{v_n}^{k_2}u_n \overline{v_n}^{k_3}\cdots.
\end{equation*}
Thus, the $(-\beta)$-expansions of real numbers which belong to gaps start by $ \sigma^i(u_{k-1})u_{k-1}u_{k-1} $ or $ \sigma^i(u_{k-1})u_ku_k$.

If $ \beta $ belongs to $ [\gamma_{n+1}, \gamma_n [ $, $ d(l_{\beta}, -\beta)$ satisfies \eqref{6}.
 As given in \eqref{eq16}
\begin{equation*}
v_k = \begin{cases}
    d_1 \cdots d_{\vert u_k \vert-1}(d_{\vert u_k \vert}-1)0  &\text{ if $ d_{\vert u_k \vert } =1 $}\\
    d_1 \cdots d_{\vert u_k \vert-2}(d_{\vert u_k \vert-1}+1) &\text{ if $ d_{\vert u_k \vert } =0 $}.
    \end{cases}
\end{equation*}
Thus, for  $ 1\leq i \leq n $, 
$ d_1 \cdots d_{2n_i-1} = u_i $ and $ d_1 \cdots d_{p_i} d_{p_i+2n_i} = v_i $ if $ i $ is odd or  $ d_1 \cdots d_{p_i} d_{p_i+2n_i}0 = v_i $ if $ i $ is even.

In fact, $ \Delta_{odd} $ or one of $ \Delta^{(i)}$ codes the support of the maximal entropy measure on $ S_{-\beta} $ endowed with the shift $ \sigma $. That is, the support is the closure of $ \Delta_{odd}^{\Z}$ or $ \Delta^{(i_0)\Z}$ for some $ i_0 $. 

\begin{rem}\label{der1}
We assume that \eqref{6} is satisfied. If there exists an integer $ i $ such that 
\begin{equation}
d_1\cdots d_{n-1}(d_n+(-1)^n)\not\in L_{\beta}, \text{ $ \forall n $, $ n\geq 2n_i-1 $}
\end{equation}
then, $(d_i)_{i \geq 1}$ is periodic with odd period. Indeed, at first, note that 
\[ d_1 \cdots d_{2n_i-1} d_1 d_2 d_3 d_4\cdots \]
is the upper (with respect to the alternating order) sequence starting by the string $ d_1 \cdots d_{2n_i-1} $. Thus, for any integer $ n $ greater than $ 2n_i$, $ d_1 \cdots d_{n-1}(d_n+(-1)^n) \not\in L_{\beta}$ means that there is no word between  $ d_1\cdots d_{2n_i-1} d_1 d_2 d_3 \cdots $ and $(d_i)_{i \geq 1}$. That is, for any $ n\geq 1$, $ d_n = d_{2n_i-1+n}$.
\end{rem}

\begin{rem}\label{der2}
Let $ \beta$ be the real number such that $ \gamma_1 < \beta \leq \gamma_0 $, and $(d_i)_{i \geq 1} $ is supposed to be non periodic with odd period. If $(d_i)_{i \geq 1}$ satisfied \eqref{6}, $ d_2 \cdots d_{2n_1-1}$ is the longest sequence of zero. It exists in $ \Delta_{odd}$ a word starting by $ d_1\cdots d_{2n_1-1}$. Indeed, from Proposition \ref{prop7},
\begin{equation*}
d(l_{\beta}, -\beta) = 1 00 (11)^{k_1}00 (1)^{k_2}00(1)^{k_3}00 \cdots. 
\end{equation*}
\begin{itemize}
\item If $ k_1\neq 0$, $ d_1\cdots d_{2n_1 - 1} = 100$ and $ 100(11)^{k_1}\in \Delta_{odd} $. 

\item If $k_1 = 0$, we can write $(d_i)_{i \geq 1 }$ on the form 
\begin{equation}
(d_i)_{i \geq 1 } = 1(00)^{t_1}1(00)^{t_2}1(00)^{t_3}1\cdots,
\end{equation}
with $t_1 \geq 2$ (we use the fact that $d(l_{\gamma_0}, -\gamma_0)\prec (d_i)_{i \geq 1}$). One has: 
\[ d_1\cdots d_{2n_1-1} = 1(00)^{t_1}.\] 
When $ t_3 \neq 0$, $ 1(00)^{t_1}1(00)^{t_2}1 $ belongs to $ \Delta_{odd}$. If $ t_3 = 0$, we have \[ 1(00)^{t_1}1(00)^{t_2}100 \in \Delta_{odd} .\] 
\end{itemize}
In both cases, there exists in $\Delta_{odd}$ a word starting by $ d_1 \cdots d_{2n_1-1} $.
\end{rem}

\begin{lem}\label{l5}
Let $ \beta $ be a real number such that $ \gamma_1 < \beta \leq \gamma_0 $. Then $ \Delta_{odd}$ codes the support of the maximal entropy measure and $ \Delta^{(i)*} \subset L_{\Delta_{odd}^*}$. 
\end{lem}

\begin{preu}
From Remark \ref{der1}, if $ (d_i)_{i \geq 1}$ is not periodic with odd period, for any integer $ k $, we can find $ n \geq 1 $ such that $ d_1 \cdots d_{n-1}(d_n+(-1)^n) \in L_{\beta} $. 

If $ n $ is even, $ d_n+1 = 1$. That is, $ d_n = 0$. From Remark \ref{der2}, there is a word $ y \in \Delta_{odd}$ starting by $ d_1\cdots d_{2n_1-1}$. The sequence $ d_1\cdots d_{n-1}y \in L_{\Delta_{odd}^*}$.
If $ n $ is odd, $ d_n-1 = 0$. That is $ d_n = 1$. Note that, $(d_i)_{i \geq 1}$ is a concatenation of 1 and $00$. The word $ d_1 \cdots d_n $ and by a string of the type $1(00)^t1$ and $(00)^t $ is not the longest sequence of zero.
\begin{equation}
d_1 \cdots d_n = d_1 \cdots d_{n-2t-2}1(00)^t1.
\end{equation}
It follows that $ d_1 \cdots d_{n-2t-2}y \in L_{\Delta_{odd}^*} $ since this word end by a word of $\Delta_{odd}$. Thus, for any admissible concatenation $ B_{k_1} \cdots B_{k_m} $, for all $ n $, $ n \geq p_{k_m}$, $ B_{k_1} \cdots B_{k_m}d_1\cdots d_n \in L_{\Delta_{odd}^*}  $. That is $ B_{k_1} \cdots B_{k_m} \in L_{\Delta_{odd}^*}  $.
\end{preu}


\begin{lem}\label{l6}
 Let $ \beta$ be a real number such that for all $ n $ in $ \N^*$, $\gamma_{n+1}\leq \beta < \gamma_n $, the support of the maximal entropy measure is coded by $ \Delta^{(n)}$.
   \end{lem}
 
 \begin{preu}
Note that  $\Delta_{odd}=\{1\}$, $\Delta^{(i)}=\{ u_i \} $ (with $ i < n $). Therefore, $ d(l_{\beta}, -\beta) \in \{u_n, v_n \}^*$ and satisfies \eqref{eq36bis}. That is 
\begin{equation*}
d(l_{\beta}, -\beta)= u_n^n v_n(u_nu_n)^{k_1}v_n(u_n)^{k_2}v_n(u_n)^{k_3}v_n \cdots.
\end{equation*} 
We obtain the same result as in Lemma \ref{l5} by changing the alphabet $\{1, 00 \} $ to $\{u_n, v_n \}$. So the language which codes the support of the maximal entropy of the $\sigma$-invariant measure contains $ u_n $. Thus, the right choice is $\Delta^{(n)}$. For $ i > n $, $\Delta^{(i)*} \subset L_{\Delta^{(n)*}}$. But this inclusion can not hold for $ i < n $. 
\end{preu}

Each set $\C $, $ \Delta_{odd}$, $\Delta^{(i)}$, $i< n$ defines a class of words forbidden in the language of the support of the maximal entropy measure. These words are: $ u_k u_k u_ku_k $ and $ u_k u_{k+1} u_{k+2} $ with $ -1 \leq k < n$ and $u_{-1}=0$. In the one side right $(-\beta)$-shift, we add $ \sigma^i(u_k)u_ku_ku_k $ and $ \sigma^i(u_k)u_{k+1}u_{k+2}$ (with $i < \vert u_k \vert $). It easy to see that one of these sequences appears in the expansion of a real taken in a gap.

\section{Lap counting function}\label{s3}
 
The lap counting function of a continuous map  $ T $ whose consists of a finite number of monotone segments (called laps) is the formal power series 
\begin{equation}
L_T(z) = \sum\limits_{ n \geq 0} L_n(T)z^n
\end{equation}
in which  $ L_0(T) = 1 $ and for all $ n \geq 1 $, $ L_n(T)$ counts the number of laps of the iterate $ T^n $. This function was been introduced by Milnor and Thurston. It is another approach to computing the zeta function.

Let $ \beta> 1$. In this section, we give lap-counting function of the $(-\beta)$- transformation and its classical properties. In the following, $ S_{-\beta}^r $ denotes the one-side right $(-\beta)$-shift ($(-\beta)$-representations of real belonging to $ I_{\beta}$). 

\begin{theo}\label{th3}
Let consider a real $ \beta> 1$ and $ T_{-\beta}$, $(-\beta)$-transformation. We set $ d(l_\beta, -\beta) = (d_i)_{i \geq 1} $. Then, the lap counting function $ L_{T_{-\beta}} $ of $T_{-\beta}$ is given by: 

\begin{equation}
L_{T_{-\beta}}(z) = \dfrac{1}{(1-z)(1-\sum\limits_{n \geq 1} (-1)^n(d^{*}_{n-1}-d^{*}_n)z^n)};
\end{equation}
where $ (d^{*}_i)_{i \geq 1} $ is defined in \eqref{D}.
\end{theo}

\begin{rem}\label{r5}
For a fixed real non integer $ \beta > 1$, the graph of $ T_{-\beta} $ consists of $ \lfloor \beta \rfloor +1 $ segments. Indeed, we can see $ I_\beta $ as union of $ \lfloor \beta \rfloor +1 $ intervals $ I_k $ defined by $ I_{0} = (-\frac{1}{\beta}+r_\beta,  r_\beta )$, $ I_k = ( -\frac{k+1}{\beta} + r_\beta, -\frac{k}{\beta}+r_\beta]$ with $ 0< k < \lfloor \beta \rfloor -1 $, and $ I_{\lfloor \beta \rfloor } = [ l_\beta, -\frac{\lfloor \beta \rfloor}{\beta} +r_\beta ]$.
\end{rem}

\begin{lem}\label{lem3}
The $(-\beta)$-transformation $ T_{-\beta}$ is affine on each interval $ I_k $.
Moreover, $ x \in I_k $ if only if $d(x, -\beta)$ starts by $ k $.
\end{lem}

\begin{preu}
Let $ x \in I_{\beta}$. 
\begin{align*}
x \in I_k & \Longleftrightarrow -\frac{k+1}{\beta} +r_\beta < x \leq -\frac{k}{\beta} +r_\beta .\\
          & \Longleftrightarrow k \leq -\beta x -l_\beta < k+1 \\
          & \Longleftrightarrow \lfloor -\beta x -l_\beta \rfloor = k.
\end{align*}
We set $ (x_i)_{i \geq 1} = d(x, -\beta) $. We know that $ x_i = \lfloor -\beta T^{i-1}_{-\beta}(x)-l_\beta \rfloor $ (see \cite{MR2534912}). So, we have proved that for all $ x \in I_k$, $ x_1 = \lfloor -\beta x - l_\beta \rfloor = k $ and $ T_{-\beta}(x) = -\beta x -k$, $ 0 \leq k \leq \lfloor \beta \rfloor $.
\end{preu} 

Note that, if $ \beta \in \N$, $ I_{\lfloor \beta \rfloor} = \{ l_\beta \} $. In this case, the graph of the $(-\beta)$- transformation  is given by $ \beta $ segments.

From the previous lemma, $ T_{-\beta}$ has $ \lfloor \beta \rfloor + 1 $ laps  if $ \beta \not\in \N$ and $ \beta $ laps otherwise.

\begin{equation}
d(-\frac{k+1}{\beta}+r_\beta, -\beta) = \cdot (k+1)d_1d_2\cdots,
\end{equation} 
with $ 0 \leq k < d_1 $ and $d(l_\beta, -\beta) = (d_i)_{i\geq 1}$. It is easy to see that $ \cdot k d_1d_2\cdots $ is the maximum (with respect to the alternating order) in the family of words of $ S_\beta $ starting by $ k $.

Throughout the rest of this section, we set $ I_{a_1a_2\cdots a_n} = \{ T_{-\beta}^{k-1}(x) \in I_{a_k}, 1 \leq k \leq n \} $. 

\begin{rem}\label{rem6}
  Let $ T $ be the map from $ S_{-\beta}$ to $ I_{\beta}$ defined by 
	\begin{equation*}
	T((x_i)_{i \geq 1}) = \sum\limits_{n \geq 1}\frac{x_n}{(-\beta)^n}.
	\end{equation*}
The map $ T $ is increasing in the sense of the alternating order (see \cite{MR2534912}). It is easy to verify that the words  $ X = x_1 \cdots x_{n-1}(x_n+1)d_1 d_2 \cdots $ and $ Y=x_1 \cdots x_{n-1}x_n0d_1d_2 \cdots $ have same image by $ T $.
\end{rem}
Indeed, 
\begin{align*}
T(X) - T(Y) &= \frac{1}{(-\beta)^n} (1+ l_\beta - r_\beta ) \\
            &= 0
\end{align*}
We have the same result if in Remark \ref{rem6}, we replace $(d_i)_{i \geq 1}$ by $ (d_{i}^*)_{i \geq 1}$.

\begin{lem}\label{lem4}
Let $ I_{a_1a_2\cdots a_n} = \{ T_{-\beta}^{k-1}(x) \in I_{a_k}, 1 \leq k \leq n \} $. It is the interval of real numbers of $ I_{\beta}$ for which the $(-\beta)$-expansions begins by the admissible word $ a_1a_2\cdots a_n $.
\end{lem}

\begin{preu}
Let $ x \in I_{a_1 a_2 \cdots a_n} $. For $ 1\leq k \leq n $, $ T_{-\beta}^{k-1}(x) \in I_{a_k} $. From the previous corollary, 
\begin{equation}
\lfloor -\beta T_{-\beta}^{k-1}(x) - l_\beta \rfloor = a_k,
\end{equation}
that is $d(x, -\beta) $ begins by $ a_1 a_2 \cdots a_n $. 

Furthermore, consider the set of $(-\beta)$-expansions of reals of $ I_{a_1 \cdots a_n }$. 
\begin{itemize}
\item[(a)] If $ a_1 \cdots a_n $ ends by a word of the type $ d_1 \cdots d_{n-1}j$ with $(-1)^n(d_n - j)<0 $ and $ j \neq d_1 $, the  $(-\beta)$-expansions of $ I_{a_1 \cdots a_n}$ endpoints are $ a_1 a_2 \cdots a_n d_1 d_2 \cdots $ and $ a_1 \cdots a_{n-1} (a_n + 1) d_1 d_2 \cdots$. In fact, 
\[ a_1 \cdots a_n 0d_1^*d_2^*d_3^* \cdots \] 
is one of the endpoints $(-\beta)$-representations starting by $ a_1 a_2 \cdots a_n $. But, this word cannot be an expansion. From Remark \ref{rem6}, both words
\[  a_1 \cdots a_{n-1} (a_n + 1) d_1 d_2 \cdots  \]
 and 
\[a_1 \cdots a_n 0d_1^*d_2^*d_3^* \cdots \]
 have the same image by $ T $. 

\item[(b)] We assume that $d_{2i} < d_1 $ for all integer $ i $ and $ a_1 \cdots a_n $ ends by a sequence of the type $ d_1 \cdots d_k $, we set 
$ a_1 \cdots a_n = a_1 \cdots a_{n-k} d_1 \cdots d_k $.
Thus, the word $ a_1 \cdots a_{n-k}d(l_\beta, -\beta) $ is a endpoint $(-\beta)$-expansion. The maximum in $ S_{-\beta} $ starting by $ d_1 \cdots d_k $ is $ d_1 \cdots d_k d(l_\beta, -\beta)$ if $ k $ is odd and  $ d_1 \cdots d_k 0 d_1^*d_2^*d_3^*\cdots $ if $ k $ is even. This last word is not a $(-\beta)$-expansion. We replace it by \[ a_1 \cdots a_{n-k}d_1\cdots d_{k-1} (d_k+1)d(l_\beta, -\beta).\]

\item[(c)] Now, we assume that $ a_1 \cdots a_n $ ends by a sequence of the type $ d_1 \cdots d_k $, $ k $ even and $ d_k = d_1$, then, the endpoints have expansions 
$a_1\cdots a_{n-k}d(l_\beta, -\beta)$ and $ a_1\cdots a_{n-k}d_1\cdots d_{k-1}d(l_\beta, -\beta).$

This last case implies that $ d(l_\beta, -\beta)$ satisfies \eqref{6}. In the language of $ S_{-\beta}$, we know that $ d_1 \cdots d_{2n_i-1} $ is always followed by $ d_1 \cdots d_{p_i}$. Thus, if $ 2n_i \leq k \leq 2n_i+p_i - 1$, the endpoints of $ I_{a_1\cdots a_n} $ have expansions $ a_1\cdots a_{n-k} d(l_\beta, -\beta) $ and $ a_1 \cdots a_{n-k}d_1\cdots d_{2n_i-1} d(l_\beta, -\beta)$. That is, for all $ t , r $ and $ s $ in $ \N $, such that $ 2n_i \leq r, s \leq 2n_i+p_i-1$, 
\begin{equation*}
I_{a_1\cdots a_t d_1 \cdots d_r} = I_{a_1\cdots a_td_1\cdots d_s }.
\end{equation*}
\end{itemize}
Suppose $ d(l_\beta, -\beta) $ periodic with odd period $ 2n-1$. There exists $ i $ such that $ d_{2n_i}= d_{2n} = d_1$ and $ d_1 \cdots d_{p_i} = 2n-1 = 2n_i-1$. (c) allows us to say that $ I_{a_1 \cdots a_kd_1\cdots d_{2n-1}}$ is reduced to the singleton $\{ T(a_1\cdots a_k) + \frac{l_\beta}{(-\beta)^k}$ \}. It is not an interval. In fact, in an infinite admissible word, $ d_1 \cdots d_{2n-1} $ is always followed by itself.
\end{preu}

\begin{rem}\label{rem7}
Let $ \beta $ be a real number strictly bigger than 1 and $ L_{\beta}$ the language of the words $ (x_i)_{i\geq 1} $ such that
\begin{equation*}
(d_i^*)_{i \geq 1} \preceq (x_{i+n})_{i \geq 1} \prec (d_{i-1}^*)_{i \geq 1},   \text{ $ \forall n \geq 0 $ and $ d_0^* = 0 $}.
\end{equation*} 
From the previous corollary, one has 
for all fixed integer $ n>0 $,
\begin{equation}
I_{\beta} = \underset{\underset{ \vert x \vert = n}{x \in L_{\beta}}}{\bigcup} I_{x}.
\end{equation}
\end{rem}
Rather than counting intervals, one can just count words of the language of $S_{-\beta}$.
This approach allows to obtain the laps of $ T_{-\beta}^n $, for all $ n \in \N^{*} $. 

\begin{preu 3}

Let $ \beta > 1$. We recall that $ d(l_\beta, -\beta) = (d_i)_{i \geq 1} $. From the formula for the factor complexity of the corrected $ (-\beta)$-shift given in Section 2 (we can also see in  \cite{NguemaNdong20161} and \cite{Ref25}), 
 and from Lemma \ref{lem3}, Lemma \ref{lem4}, and Remark \ref{rem7}, one has
\begin{equation*}
L_n = \sum\limits_{i = 1}^{n}(-1)^i(d^*_{i-1}-d^*_i)L_{n-i} + 1.
\end{equation*}
In fact, $ L_n $ counts the number of words of length $ n $ in the language of $ \tilde{S}_{-\beta}$. Then, by simple calculus in the open unit disk except in $ \frac{1}{\beta}$, 
\begin{equation*}
L_{T_{-\beta}}(z) = \left( \sum\limits_{n \geq 1}(-1)^n(d_{n-1}^*-d_n^*)z^n\right) L_{T_{-\beta}}(z) + \sum\limits_{n \geq 0}z^n.
\end{equation*}
Hence the result follows. 
\end{preu 3}

\begin{exple}
 The following figures represent the graphs of $ T_{-2.5}$, $T_{-2.5}^2 $ and $ T_{-2.5}^3$ respectively in $ I_{2.5} = [-\frac{2.5}{3.5}, \frac{1}{3.5} ) $. 
 
\begin{figure}[!ht]
\centering
\includegraphics[width=10cm, height=3cm]{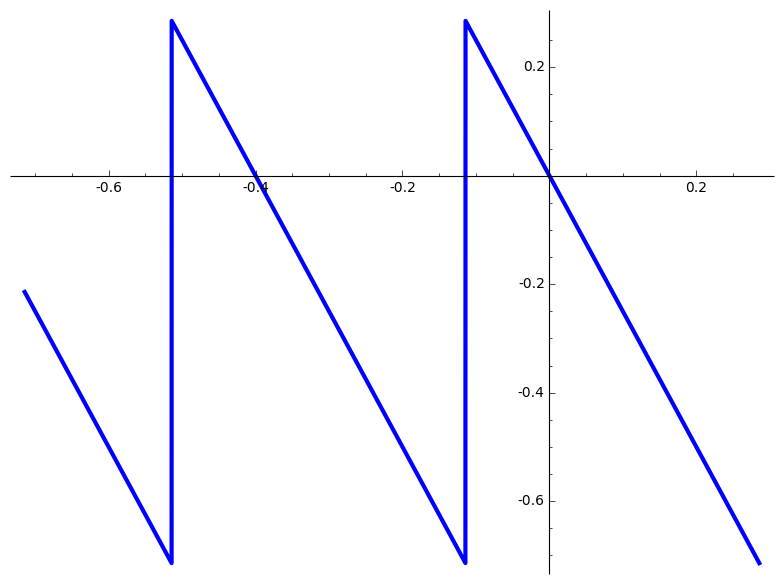}
\caption{$T_{-\beta}$ for $ \beta = 2.5 $}
\end{figure}

\begin{figure}[!ht]
\centering
\includegraphics[width=10cm, height=3cm]{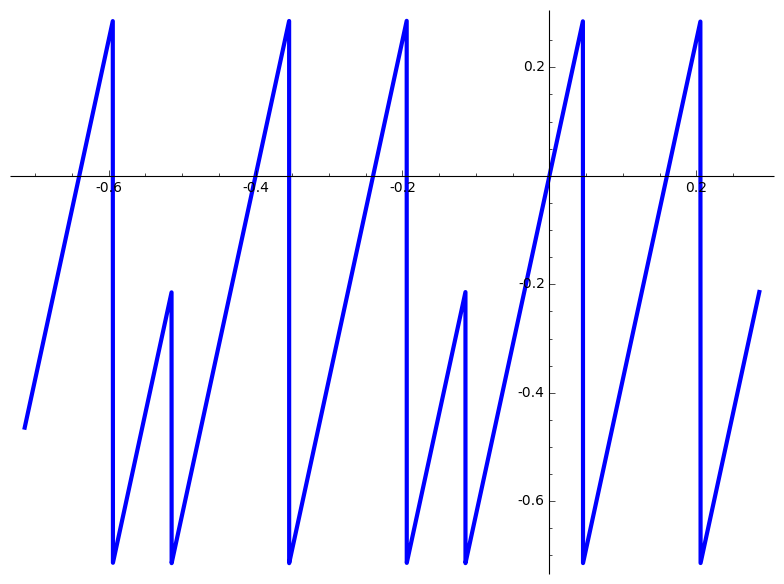}
\caption{ $T_{-\beta}^2$ for $ \beta = 2.5 $}
\end{figure}

\begin{figure}[!ht]
\centering
\includegraphics[width=10cm, height=3cm]{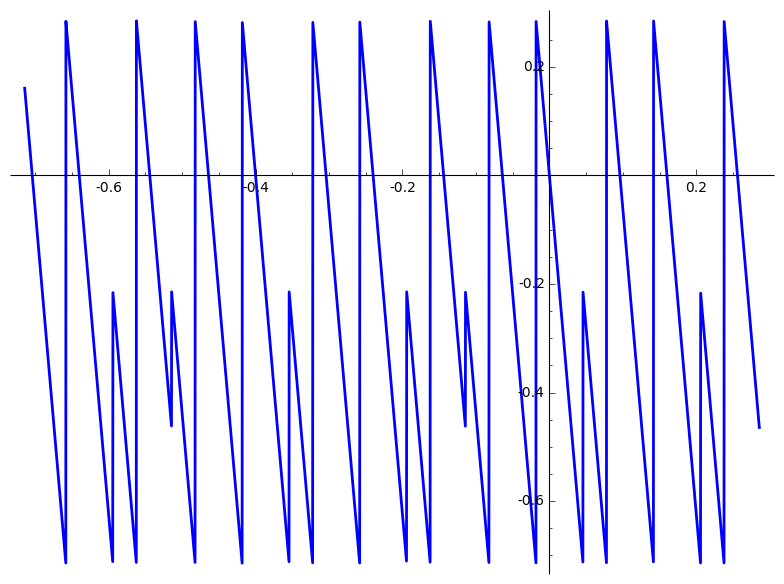}
\caption{ $T_{-\beta}^3$ for $ \beta = 2.5 $}
\end{figure}

 We set $d(l_\beta, -\beta) = (d_i)_{i \geq 1}$ and we verify easily that $ d_1 = 2$, $d_2 = 1$ and $ d_3 = 1$. The laps of $ T_{-\beta}^i $ correspond to the number of oblique segments in the different graphs. These oblique segments allow to determine the number of words of length $ i $ in the $(-2.5)$-shift language. 
\begin{itemize}
\item For $ i=1 $, we have three laps (oblique segments) and three admissible words of length 1: 0, 1 and 2 ($ d_1 = 2$). 
\item If $ i = 2$, $ d_1 d_2 = 21$; there are 8 oblique lines and then 8 admissible words of length 2: 21, 22, 10, 11, 12, 00, 01, 02.
\item When $ i = 3$, $ d_1 d_2 d_3 = 211$; we count twenty laps, then there exist 20 admissible words of length $3 $: 211, 210, 222, 221, 102, 101, 100, 112, 111, 110, 122, 121, 002, 001, 000, 012, 011, 010, 022, 021.
\end{itemize}  
 \end{exple}

\section{Zeta function}

The notion of dynamical zeta function was been introduced by M. Artin and B. Mazur in 1965. We consider a diffeomorphism $ \delta $ on a compact space such that all of its iterates $ \delta^n $ have isolated fixed points. The zeta function associated to $ \delta $ is given by: 
\begin{equation}
 \zeta_{\delta}(z) = exp( \sum_{k \geq 1} \frac{\sharp Fix(\delta^k)}{k}z^k)
 \label{(27)}
\end{equation}
where $ \sharp Fix (\delta^k) $ counts the number of fixed points of $ \delta^k $, by analogy with the geometric zeta function.

In 1994, Lepold Flatto, Jeffrey Lagarias and Bjorn Poonen (see \cite{MR1279470}) dealt with the zeta function of the $ \beta $-transformation. 
They consider the application from $ [0, 1) $ to $ [0, 1) $ defined by: 
\begin{equation*}
 T_{\beta} : x \mapsto \left\lbrace \beta x \right\rbrace  \text{ for $ \beta > 1 $ }
\end{equation*}
where $ \left\lbrace x \right\rbrace $ denotes the fractional part of $ x $.
The associated zeta function is: 
\begin{equation}
 \zeta_{\beta}(z) = exp ( \sum_{k \geq 1} \frac{p_k}{k}z^k )
 \label{(28)}
\end{equation}
where $ p_k $ counts the number of fixed points of $ T_{\beta}^k $. In other words, $ p_k $ is the number of periodic admissible sequences $ .x_1x_2\cdots $ with period $ k $. We denote by $ d_{\beta}(x) $ the expansion of $ x $ in base $ \beta $. 

\begin{equation*}
  x_1x_2\cdots = d_{\beta}(x), T_{\beta}^k (x) = x \Rightarrow d_{\beta}(x)=d_{\beta}(T_{\beta}^k(x)).
\end{equation*}
And then, $ ( x_{k+i} )_{i \geq 1} = (x_i)_{i \geq 1} $ since the expansion is unique for each number and all sub-words of an admissible word is an expansion. 

After the introduction of the $(-\beta)$-expansion in 2009 by Ito and Sadahiro in \cite{MR2534912}, in the following sentences, we focus our study on the determination of the zeta functions of the $(-\beta)$-transformation and the one of the  $(-\beta)$-shift endowed with the shift. 
\vspace{1cm}

\subsection[Zeta function and exhaustive code]{Zeta function of symbolic dynamical system defined by an exhaustive prefix code}
\vspace{0.5cm}

Let $ X $ be a symbolic dynamical system. Suppose $ X $ coded by $ C $. 

\begin{defi}
The code $ C $ is said \textit{ exhaustive } if all periodic word $ P $ can be written uniquely as: 

\begin{equation*}
  P = a_1a_2\cdots a_s x_{1,1}x_{1,2}\cdots x_{1,k_1}x_{2,1} \cdots x_{2,k_2} \cdots x_{h-1,1}\cdots x_{h-1,k_{h-1}}b_1 \cdots b_r
\end{equation*}
where
\begin{align*} x_1 & = x_{1,1} \cdots x_{1,k_1} &&\in C \\
                 x_2 & = x_{2,1} \cdots x_{2,k_2} &&\in C \\
               & \vdots \\
               x_{h-1} & = x_{h-1,1} \cdots x_{h-1, k_{h-1} } &&\in C \\
                   x_h & = b_1 b_2 \cdots b_r a_1a_2 \cdots a_s &&\in C.
\end{align*}
\end{defi}
$ P $ and $ x_1 x_2 \cdots x_h $ have the same orbit. 

\begin{theo}\label{t3}
 Let $ X $ be a coded system defined by an exhaustive prefix code $ C $. 
Then, if $ p_n $ counts the number of periodic words of period $n$ in $ X $ the associated zeta function is defined by 
\begin{equation}
\begin{aligned}
\zeta_X (t) &= exp( \sum_{n \geq 1} \frac{p_n}{n}t^n ) \\
&=\frac{1}{1 - \underset{n \geq 1}{\sum}b_n t^n} 
\label{(29)}
\end{aligned}
\end{equation}
where $ b_n $ counts the number of words of length n in $ C $.
\end{theo}

\begin{preu}
 Observe that, in the meaning of formal power series,
\begin{equation*}
\begin{aligned}
\log \dfrac{1}{1-\sum\limits_{n \geq 1} b_n t^n} &= -\log \left(1-\sum\limits_{n\geq 1}b_nt^n\right) \\
&=\underset{n \geq 1}{\sum}b_nt^n + \frac{1}{2}\left( \underset{n \geq 1}{\sum} b_nt^n\right)^2 + \frac{1}{3}\left(\underset{n \geq 1}{\sum}b_nt^n\right)^3 +\cdots\\
&=\sum\limits_{k\geq 1}\frac{1}{k}\left( \sum\limits_{n \geq 1}b_nt^n\right)^k. 
\end{aligned}
\end{equation*}    

Let $ (\delta_{n,k})_{n\geq k \geq 1} $ be the sequence of positive integers such that $ \delta_{n,k} $ counts the number of periodic words of length $ n $ having the same orbit than a product of $ k $ pieces of the exhaustive prefix code. Thus, $ p_n = \sum_{k=1}^n \delta_{n,k} $. In the sense of formal power series, we rewrite $ \sum_{n \geq 1}\frac{p_n}{n}t^n $ as: 

\begin{equation*}
  \sum_{n \geq 1}\frac{p_n}{n}t^n = \sum_{n \geq 1}\frac{\delta_{n,1}}{n}t^n+\cdots + \sum_{n \geq k}\frac{\delta_{n,k}}{n}t^n+ \cdots . 
\end{equation*}
Let $ x_1, x_2, \cdots, x_k $ k pieces of $ C $ with $ \vert x_1x_2\cdots x_k \vert = n $ and $ p, h $ two integers such that $ x_1 x_2 \cdots x_k = (x_1x_2\cdots x_p)^h $ where $ p $ 
is minimal, that is $ x_1 \cdots x_p $ denotes the smallest word (in size) with orbit $ \overline{x_1x_2 \cdots x_k}$.

\begin{equation*}
  k = hp \text{ and } \vert x_1x_2\cdots x_k \vert = \sum_{i = 1}^{k} \vert x_i \vert = h \sum_{i = 1}^{p} \vert x_i \vert. 
\end{equation*}
Then,

\begin{align*}
\delta_{n,k} & = \underset{\overset{(i_1, i_2, \cdots i_p)}{\underset{ph=k}
{h\sum_{j=1}^pi_j=n}}}{\sum}(\sum_{j=1}^{p}i_j)(b_{i_1}\cdots b_{i_p})^h \\
                                 & =  \underset{\overset{(i_1, i_2, \cdots i_p)}
{\underset{ph=k}{h\sum_{j=1}^pi_j=n}}}{\sum}\frac{n}{h}(b_{i_1}\cdots b_{i_p})^h. 
\end{align*}
Thus, 

\begin{align*} 
 \underset{n \geq k} {\sum} \frac{1}{n}\delta_{n,k}t^n & = \underset{n \geq k}{\sum}
\underset{\overset{(i_1, i_2, \cdots i_p)}{\underset{ph=k}{h\sum_{j=1}^pi_j=n}}}
{\sum}\frac{1}{h}(b_{i_1}\cdots b_{i_p})^h t^n \\
                                                       & =\frac{1}{k}\underset{n \geq k}
{\sum}\underset{\overset{(i_1, i_2, \cdots i_p)}{\underset{ph=k}{h\sum_{j=1}^pi_j=n}}}
{\sum}\frac{k}{h}(b_{i_1}\cdots b_{i_p})^h t^n \\
                                                       & =\frac{1}{k}\underset{n \geq k}
{\sum}\underset{\overset{(i_1, i_2, \cdots i_p)}{\underset{ph=k}{h\sum_{j=1}^pi_j=n}}}{\sum}p(b_{i_1}\cdots b_{i_p})^h t^n
\end{align*}
$ p $ minimal,  $ ph = k\text{ and } h \sum_{j=1}^p i_j = n$, the integer $ p (b_{i_1}\cdots b_{i_p})^h $ counts the periodic words resulting to the circular permutations of concatenations of $ k $ monotonic piecewises of the code. Then, regardless of the commutative property of the multiplication in $ \R $, and in the sense of formal power series, 
\begin{equation*}
  \underset{n \geq k} {\sum} \frac{1}{n}\delta_{n,k}t^n = \frac{1}{k} \underset{n \geq k}
{\sum}{\underset{ \overset{k}{\underset{j=1}{\sum}}i_j = n}{\sum}}b_{i_1}\cdots b_{i_k}t^n, 
\end{equation*}
but
\begin{equation*}
 \underset{n \geq k}{\sum}{\underset{ \overset{k}{\underset{i=1}{\sum}}i_j = n}{\sum}}b_{i_1}\cdots b_{i_k}t^n
  =  (\underset{ n \geq 1}{\sum}b_nt^n)^k .  
\end{equation*}
That is,
\begin{equation*}
  \underset{n \geq k}{\sum} \frac{1}{n}\delta_{n,k}t^n =\frac{1}{k}(\underset{ n \geq 1}{\sum}b_nt^n)^k. 
\end{equation*}
Hence,
\begin{equation*}
 \underset{n \geq 1} {\sum}\frac{p_n}{n}t^n = \underset{n \geq 1}{\sum}b_nt^n + \frac{1}{2} \left(\underset{n \geq 1}{\sum} b_nt^n\right)^2 + \frac{1}{3}\left(\underset{n \geq 1}{\sum}b_nt^n\right)^3 +\cdots                                           
\end{equation*}
                             
\end{preu}

Theorem \ref{t3} reveals an important property of coded systems: the density of the set of periodic points.

\begin{exple}
 Let $ \beta $ be a real number, $ \beta > 1 $. Let $ X_\beta $ be the $ \beta$-shift and $ (a_i)_{i \geq 1} $ the expansion of 1 in base $ \beta$. We assume that $ \beta $ is not a simple $ \beta$-number and we set $ C_\beta = \left\lbrace a_1\cdots a_{k}i, k\in \N, 0\leq i \leq a_{k+1}-1 \right\rbrace $.

The $ \beta$-shift $ X_\beta $ is coded by $ C_\beta $ which is an exhaustive prefix code. The integer $ a_k $ counts the number of pieces of length $k$ in $ C_{\beta}$. Thus, the zeta function associated to $ X_\beta $ is given by: 
\begin{equation*}
 \zeta_{X_{\beta}}(z) = \dfrac{1}{1-\sum\limits_{n \geq 1}a_nz^n}. 
\end{equation*}

\end{exple}

\subsection{Zeta function of the negative beta-transformation}
\vspace{0.5cm}

We consider a real number $ \beta > 1$. Recall that the $ (-\beta)$-transformation $ T_{-\beta}$ denotes the map from $ I_{\beta} = [\frac{-\beta}{\beta+1}, \frac{1}{\beta+1} ) $ into itself defined by:
\begin{equation*}
T_{-\beta}(x) = -\beta x - \lfloor -\beta x + \dfrac{\beta}{\beta+1} \rfloor.
\end{equation*}
The aim of this section is to determine the zeta function $ \zeta_{-\beta}$ of the map $ T_{-\beta}$. We know that each number has a $(-\beta)$-representation in $\tilde{S}_{-\beta}$. Moreover, it is easy to see that a real $ x $ is a fixed point  of $ T_{-\beta}^k$ if only if the $(-\beta)$-representation of $ x $ in $\tilde{S}_{-\beta}$ is periodic with period dividing $ k $.

\begin{theo}\label{t2}
Let $ \beta $ be a real number strictly greater than 1, $ (d_i)_{i \geq 1}$ the $(-\beta)$-expansion of $ l_\beta = -\frac{\beta}{\beta + 1}$ and $ \zeta_{-\beta} $ the Zeta function of the $(-\beta)$-transformation. Then, in the ball of radius $ \frac{1}{\beta}$ and center 0, 
\begin{itemize}
\item if $(d_i)_{i \geq 1}$ is not periodic 
\begin{equation*}
\zeta_{-\beta}(z) = \dfrac{1+z}{1-\sum\limits_{n \geq 1}(-1)^n(d_{n-1}-d_{n})z^n};
\end{equation*}
\item if $(d_i)_{i \geq 1}$ is periodic with period $ k $,  
\begin{equation*}
\zeta_{-\beta}(z) = \dfrac{1+z}{(1-z^k)(1-\sum\limits_{n \geq 1}(d_{n-1}^* - d_n^*)(-z)^n)}.
\end{equation*}
\end{itemize}
\end{theo}

If we consider a real $\beta> 1$ and $(a_i)_{i \geq 1}$ the $\beta$-expansion of 1, the zeta function of the $\beta$-shift, determined by Leopold Flatto, Jeffrey Lagarias and Bjorn Poonen in \cite{MR1279470}, is given by
\begin{equation*}
\zeta_{\beta}(z) = \frac{1}{1-\sum\limits_{n \geq 1} a_nz^n}.
\end{equation*}
So, we remark some similarities between this zeta function and that of the $(-\beta)$-shift given in the previous theorem. For instance, $ \frac{1}{\beta} $ is a pole of these both functions.
\vspace{0.5cm}

\begin{lem}\label{lem2}
Let $\beta> 1$ and $ d(l_\beta, -\beta) = (d_i)_{i \geq 1}$ periodic with period $ h $.
\begin{equation*}
d(l_\beta, -\beta) = \cdot \overline{ d_1\cdots d_h }.
\end{equation*}
Then, $ d_h \neq 0 $.
\end{lem}
\begin{preu}
Suppose $ d_h = 0$. Since $ (d_{h+i})_{i \geq 1} = (d_i)_{i \geq 1} $
\begin{equation*}
T_{-\beta}^n(l_\beta) = \sum\limits_{n  \geq 1}\frac{d_{n+i}}{(-\beta)^i}\Rightarrow \sum\limits_{i \geq 1}\frac{d_{h+i}}{(-\beta)^i} = l_\beta.
\end{equation*}
Furthermore, 
\begin{equation*}
d_{h-1} = \lfloor -\beta T_{-\beta}^{h-2}(l_\beta)-l_\beta \rfloor .
\end{equation*}
Since $ T_{-\beta}^{h-2}(l_\beta) = \sum\limits_{i \geq 1}\frac{d_{h-2+i}}{(-\beta)^i}$ and $ d_h = 0$, it follows that
\begin{equation*}
d_{h-1} = \lfloor d_{h-1} - \frac{1}{\beta} \sum\limits_{i \geq 1}\frac{d_{h+i}}{(-\beta)^i} -l_\beta \rfloor.
\end{equation*}
Thus, $ d_{h-1} = d_{h-1} + 1$. This is absurd. Then, $ d_h\neq 0$.
\end{preu}

As consequence, $ (d_{i-1})_{i \geq 1}$ (with $ d_0 = 0$ ) is not periodic. Then, in addition to periodic words listed in the previous paragraphs, if $ d(l_\beta, -\beta) $ is periodic with even period, we should take account of circular permutations of $ (d_1d_2 \cdots d_{2p})^k$. For each integer $ k $, there are $ 2p $  words.

Let $ X_{-\beta}^r$ be the set of $(-\beta)$-expansions of real numbers which belongs to $ I_{\beta}$. We know that the correspondence $ X_{-\beta}^r \leftrightarrow I_{\beta}$ is one to one. In fact, each real number has one and only one $(-\beta)$-expansion. Moreover, $ X_{-\beta}^r$ is invariant by the shift (all infinite sub-words of a $(-\beta)$-expansion is a $(-\beta)$-expansion). Let $ x \in I_{\beta} $ and $d(x, -\beta) = (x_i)_{i \geq 1}$. Then, $ d(T_{-\beta}^n(x), -\beta) = (x_{i+n})_{i\geq 1}$. So, if $ x $ is a fixed point of $ T_{-\beta}^{n}$, then $ d(x, -\beta)$ is periodic with period dividing $ n $. The number of fixed points of $T_{-\beta}^n $ equals the number of periodic orbits in $X_{-\beta}^r$ with period dividing $ n $. 

Note that in $ S_{-\beta}^d$, the sequences which are not $(-\beta)$-expansions end by $(d_{i-1}^*)_{i \geq 1}$ (with $ d_0^* = 0 $). 

When $(d_i)_{i \geq 1}$ is not periodic with odd period $(d_{i-1}^*)_{i \geq 1} = (d_{i-1})_{i \geq 1}$ is not periodic. In this case, the periodic orbits of $S_{-\beta}^r$ belong to $ X_{-\beta}^r$. 

When $(d_i)_{i \geq 1}$ is periodic with odd period, $(d_{i-1}^*)_{i \geq 1}$ is periodic. The periodic orbits of $S_{-\beta}^r$ belong to $ X_{-\beta}^r$ except $(d_{i-1}^*)_{i \geq 1}$. 

\vspace{2mm}

\begin{rem}\label{rmflo}
In Section 2, we distinguished in $ \tilde{S}_{-\beta} $ three types of admissible concatenations and then three types of periodic words: 
\begin{itemize}
\item at first, there are the concatenations in $ \C $;
\item secondly, the concatenations of words of $ \Delta_{odd}$; 
\item and finally, there are products of words of $ \Delta^{(i)} $, for all positive integer $ i $. 
\end{itemize}
In other hands, the periodic orbits of $ \tilde{S}_{-\beta}$ are circular permutations of sequences of $ \C^{\Z} \cup \Delta_{odd}^{\Z} \cup \left( \underset{i \geq 1}{\bigcup} \Delta^{(i)\Z} \right) $, if $ (d_i^*)_{i \geq 1} $ not periodic. When $(d_i)_{i\geq 1}$ is periodic with period $ 2p$, we add the circular permutations of  $(d_i)_{i\geq 1}$. If the period of $(d_i)_{i\geq 1}$ is $2p - 1$, $(d_i)_{i\geq 1}$ does not belong to $\tilde{S}_{-\beta}$. Let $ \mathbb{P}(X) $ be the set of periodic orbit of $ X $. One has:

\begin{equation}
\mathbb{P}(\tilde{S}_{-\beta}) =
 \begin{cases}
\C^{\Z} \cup \Delta_{odd}^{\Z} \cup \left( \underset{i \geq 1}{\bigcup} \Delta^{(i)\Z} \right) &\text{ if $(d_i^*)_{i \geq 1}$ is not periodic}\\
\C^{\Z} \cup \Delta_{odd}^{\Z} \cup \left( \underset{i \geq 1}{\bigcup} \Delta^{(i)\Z} \right) \cup\{(d_i^*)_{i \geq 1}\} &\text{ if $(d_i^*)_{i \geq 1} $ is periodic }\label{PS}
\end{cases}
\end{equation}
\begin{equation}
\mathbb{P}(X_{-\beta}^r) =
 \begin{cases}
\C^{\Z} \cup \Delta_{odd}^{\Z} \cup \left( \underset{i \geq 1}{\bigcup} \Delta^{(i)\Z} \right) &\text{ if $(d_i)_{i \geq 1}$ is not periodic}\\
\C^{\Z} \cup \Delta_{odd}^{\Z} \cup \left( \underset{i \geq 1}{\bigcup} \Delta^{(i)\Z} \right) \cup\{(d_i)_{i \geq 1}\} &\text{ if $(d_i)_{i \geq 1} $ is periodic}\label{PT}
\end{cases}
\end{equation}
\end{rem}
Let $ p_n $  be number of fixed points of $\sigma^n$ in $\tilde{S}_{-\beta} $ and $ q_n $ the number of fixed points of $ T_{-\beta}^n$, according to the previous remark, $ p_n = q_n $ if $(d_i)_{i \geq 1} $ is not periodic with odd period. If $(d_i)_{i \geq 1} $ is periodic with odd period
\begin{equation}
q_n = \begin{cases}
        p_n, &\text{ if $ 2p-1 \not\vert n$, $2p\not\vert n $}, \\
        p_n + 2p-1, &\text{ if $ 2p-1\vert n$, $2p\not\vert n $},\\
        p_n-2p, &\text{ if $ 2p-1 \not\vert n$, $2p\vert n $}, \\
        p_n-1, &\text{ if $ 2p-1 \vert n$, $2p\vert n $}.\label{qn}
        \end{cases}
\end{equation}

According to \eqref{qn}, when $(d_i)_{_ \geq 1}$ is periodic with odd period, we have the following relation between the zeta function of $ T_{-\beta}$ (denoted by $\zeta_{-\beta}$) and that of $ \tilde{S}_{-\beta}$ (denoted by $ \zeta_{\tilde{S}_{-\beta}} $):
\begin{equation}
(1-t^{2p-1})\zeta_{-\beta}(t) = (1-t^{2p})\zeta_{\tilde{S}_{-\beta}}(t)\label{Rzszt}.
\end{equation}
We use the fact that
\begin{equation}
\sum\limits_{n \geq 1}\frac{q_n}{n}t^n = \sum\limits_{n \geq 1}\frac{p_n}{n}t^n - \sum\limits_{\underset{\text{$2p-1\not\vert k$}}{k\geq 1}}\frac{1}{k}t^{2pk}+\sum\limits_{\underset{\text{$2p \not| k$}}{k\geq 1}}\frac{1}{k}t^{(2p-1)k}-\sum\limits_{k\geq 1}\frac{1}{2p(2p-1)k}t^{2p(2p-1)k}
\end{equation}
and
\begin{equation}
\sum\limits_{k\geq 1}\frac{1}{2p(2p-1)k}t^{2p(2p-1)k} =  \sum\limits_{\underset{2p-1\vert k}{k\geq 1}}\frac{1}{k}t^{2pk}-\sum\limits_{\underset{2p\vert k}{k\geq 1}}\frac{1}{k}t^{(2p-1)k}.
\end{equation}

\begin{preu 1}
The sets $ \C^{\Z} $, $ \Delta_{odd}^{\Z} $ and $ \Delta^{(i)\Z} $ are coded by $ \C $,  $ \Delta_{odd} $ and $ \Delta^{(i)} $ respectively. 
Let $q_n$, $ p_{\C, n}$, $p_{0,n}$ and $ p_{i, n}$ count the number of fixed points $ T_{\beta}^n $, $\sigma^n $ in $ \C^{\Z} $, $ \Delta_{odd}^{\Z} $ and $ \Delta^{(i)\Z} $ respectively. 
If $ (d_i^*)_{i \geq 1} $ is not periodic, 
\begin{equation}
q_n = p_{\C, n} + p_{0, n}+\sum\limits_{i \geq 1}p_{i, n};
\end{equation}
if $ (d_i)_{i \geq 1} $ is periodic with period $ 2p $, 
\begin{equation}
q_n = \begin{cases} 
         p_{\C, n} + p_{0, n}+\sum\limits_{i \geq 1}p_{i, n} &\text{ if $ 2p \not\vert n$ } \\
          p_{\C, n} + p_{0, n}+\sum\limits_{i \geq 1}p_{i, n} + 2p &\text{ if $ 2p \vert n$ }.
       \end{cases}
\end{equation}
If  $ (d_i)_{i \geq 1} $ is periodic with period odd $ 2p-1 $, 
\begin{equation}
q_n = \begin{cases} 
         p_{\C, n} + p_{0, n}+\sum\limits_{i \geq 1}p_{i, n} &\text{ if $ 2p-1 \not\vert n$ } \\
       p_{\C, n} + p_{0, n}+\sum\limits_{i \geq 1}p_{i, n} + 2p-1 &\text{ if $ 2p-1 \vert n$ }.
       \end{cases}
\end{equation}
Thus, $ \zeta_{-\beta}$ is product of elementary Zeta functions of  $ \C^{\Z} $, $ \Delta_{odd}^{\Z} $ , $ \Delta^{(i)\Z} $ and also $\{ d_1 \cdots d_{k} \}^{\Z}$ (if $ (d_i)_{i \geq 1} $ is periodic with period $ k$ ). From Theorem \ref{t3}, 
\begin{equation}
\zeta_{-\beta}(z) = \frac{1}{1-\sum\limits_{x \in \C} z^{\vert x \vert} } \frac{1}{1-\sum\limits_{x \in \Delta_{odd}} z^{\vert x \vert}}\left(  \prod\limits_{i \geq 1} \frac{1}{1-\sum\limits_{x \in \Delta^{(i)}} z^{\vert x \vert}}\right) \label{ef21}
\end{equation}
if $(d_i)_{i \geq 1} $ is not periodic; or 
\begin{equation}
\zeta_{-\beta}(z) = \frac{1}{1-\sum\limits_{x \in \C} z^{\vert x \vert} } \frac{1}{1-\sum\limits_{x \in \Delta_{odd}} z^{\vert x \vert}}\left( \prod\limits_{i \geq 1} \frac{1}{1-\sum\limits_{x \in \Delta^{(i)}} z^{\vert x \vert}}\right)\frac{1}{1-z^{k}}\label{ef22}
\end{equation}
if $(d_i)_{i \geq 1} $ is periodic with period $ k $. We conclude thanks to \eqref{eq39}.

\end{preu 1}

From \eqref{PS}, if $ (d_i)_{i \geq 1}$ is periodic with odd period $ 2p-1$,
\begin{equation}
\zeta_{\tilde{S}_{-\beta}}(z) = \frac{1}{1-\sum\limits_{x \in \C} z^{\vert x \vert} } \frac{1}{1-\sum\limits_{x \in \Delta_{odd}} z^{\vert x \vert}}\left( \prod\limits_{i \geq 1} \frac{1}{1-\sum\limits_{x \in \Delta^{(i)}} z^{\vert x \vert}}\right)\frac{1}{1-z^{2p}}.\label{zs}
\end{equation}
With \eqref{zs} and \eqref{ef22}, we find again the relation given in \eqref{Rzszt}.

The previous theorem can be proved also just using the Lap-counting function. Indeed, from the equation \eqref{eq39} and the formula of the factor complexity, the coefficients of the power series expansions of $ \frac{1}{1-\sum\limits_{x \in \C} z^{\vert x \vert} } $, $ \frac{1}{1-\sum\limits_{x \in \Delta_{odd}} z^{\vert x \vert}}$ and $ \frac{1}{1-\sum\limits_{x \in \Delta^{(i)}} z^{\vert x \vert}}$ for all $ i $, count the fixed points of iterates of $ T_{-\beta}$ except the orbit of the left end point of $ I_{\beta}$ when its expansion is periodic.

Furthermore, \ref{ef21} and \ref{ef22} provide us an interesting information on the influence of gaps in the interval on the $(-\beta)$-transformation zeta function for $ \beta < \frac{1+\sqrt{5}}{2}$. They correspond to factors in the denominator of the zeta function.
 
\begin{rem}\label{r7}

From Section \ref{s3} and Theorem \ref{t2}, we have the following relation between the zeta-function of the $(-\beta)$-transformation and its lap-counting function: 
\begin{align*}
       \zeta_{-\beta}(z) &= (1-z^2)L_{T_{-\beta}} &&\text{ if $ d(l_\beta, -\beta)$ is non-periodic } \\
(1-z^k)\zeta_{-\beta}(z) &= (1-z^2)L_{T_{-\beta}} &&\text{ if $d(l_\beta, -\beta)$ is periodic with period $ k $. }
\end{align*}
\end{rem}

\begin{exple}
Let $ \gamma_0 $ be the golden ratio: $\gamma_0 = \frac{1+\sqrt{5}}{2} $, $ d(l_{\gamma_0}, -\gamma_0) = 1 \overline{0} $. The zeta function is given by two types of periodic words: the sequences of zero and words $ x_1 x_2 \cdots x_n $ such that $ x_i \in \{10^{2n}, n \in \N \} $. $ \vert 0 \vert = 1$, $ \vert 10^{2n} \vert = 2n+1$. Moreover, the $(-\gamma_0)$-shift and  the $(-\gamma_0)$-transformation have the same zeta function, since $ d(l_{\gamma_0}, -\gamma_0)$ is not periodic. Thus
\begin{align*}
\zeta_{-\gamma_0}(z) &= \dfrac{1}{(1-z^{\vert 0 \vert})(1-\sum\limits_{n \geq 0} z^{\vert 10^{2n} \vert })} \\
                  &= \dfrac{1+z}{1-z-z^2}.
\end{align*}
Moreover, $ \sum\limits_{n \geq 1} (-1)^n(d_{n-1}-d_n)z^n = z+z^2 $. Then, 
\begin{equation*}
L_{T_{-\gamma_0}}(z) = \dfrac{1}{(1-z)(1-z-z^2)}.
\end{equation*}
We verify easily that $  \zeta_{-\gamma_0} = (1-z^2)L_{T_{-\gamma_0}}$. 
\end{exple}

The previous remark allows to observe that the relation between the zeta function of $T_{-\beta}$ and its lap-counting function differs a little bit from that of the $\beta$-transformation and the associated lap-counting function. Indeed, denote by $ \zeta_{\beta} $ and $ L_{T_{\beta}} $ the zeta function and the lap-counting function of the $\beta$-transformation, according to \cite{MR1279470}, $ \zeta_\beta = (1-z)L_{T_{\beta}} $ if $ \beta $ is not simple $\beta$-number.

\subsection{Zeta-function of the negative beta-shift}
\vspace{0.5cm}

Consider a real number $ \beta > 1$. We have seen that when $ d(l_{\beta}, -\beta)$ is not periodic with odd period, $(d_{i-1}^*)_{i \geq 1}$ is not periodic and for a periodic orbit $ \overline{x_1x_2\cdots x_n} $, 
\begin{equation*}
 d_1d_2\cdots \prec \overline{x_kx_{k+1} \cdots x_nx_1\cdots x_{k-1}} \prec 0d_1^{*}d_2^{*}\cdots \text{ $ \forall k, 1 \leq k \leq n $. } 
\end{equation*}
In other words, all periodic word is a $(-\beta)$-expansion. And then, $ S_{-\beta}$ and $ T_{-\beta}$ have the same zeta function.

However, if $ (d_i)_{i \geq 1}$ is periodic with odd period $ p $, $ 0d_1^*d_2 ^* \cdots $ is periodic too (with period $ p+1$). But, it is not an expansion in base $ -\beta $. Let $ p_n $ counts the number of fixed points of $ T_{-\beta}^n$ and $ \tilde{p}_n $ the number of periodic words with period dividing $ n $ in $ S_{-\beta}$. We have the following result
\begin{equation*}
\tilde{p}_n = \begin{cases}
                           p_n &\text{ if $ p+1 \not\vert n $} \\ 
                           p_n + (p+1) & \text{ if $ p+1 \vert n $}.
                                     \end{cases} 
\end{equation*}
$ p+1 $ counts the circular permutations of the sequences $ \overline{d_1^*\cdots d_{p-1}(d_p-1)0}$. We denote by $ \tilde {\zeta}_{-\beta}$ the zeta function of the $ (-\beta)$-shift. Then, 
\begin{equation}
(1-z^{p+1})  \tilde {\zeta}_{-\beta}(z) = \zeta_{-\beta}(z).
\end{equation}
In short, $ \forall z \in B(0, \frac{1}{\beta})$,
\begin{equation*}
\zeta_{-\beta}(z) = \begin{cases} 
                             (1-z^{p+1})\tilde{\zeta}_{-\beta}(z) & \text{ if $ d(l_\beta,-\beta)=\overline{d_1 \cdots d_p} $, $p $ odd } \\
                             \tilde{\zeta}_{-\beta}(z) &\text{ otherwise}.
                             \end{cases}
\end{equation*}

\section*{Conclusion}

Finally, we have seen that for a real number $ \beta> 1$, the associated $(-\beta)$-shift is coded if only if $ \beta \geq \frac{1+\sqrt{5}}{2} $ and the $(-\beta)$-expansion of $ l_{\beta}=- \frac{\beta}{\beta+1} $ is not periodic with odd period. The non-coded case is due to the existence of intransitive words in the language of the system. In the periodic case with odd period, shall we say $ d(l_\beta, -\beta) = \overline{d_1 \cdots d_{2p-1}}$, the word $ d_1 \cdots d_{2p-1} $ is intransitive. For $ \gamma_n \geq \beta $, $ u_k u_ku_ku_k $ is an intransitive word, with $ k < n-1 $, where $d(l_{\gamma_n}, -\gamma_n) = u_n \overline{u_{n-1}}$, $ u_k = \phi^k(1)$, $ \phi(1) = 100$, and $ \phi(0) = 1 $. 

Moreover, the zeta-functions of the $(-\beta)$- shift and $\beta$- shift (determined by Leopold Flatto, Jeffrey C. Lagarias and Bjorn Poonen in \cite{MR1279470}) have some similarities: $ \frac{1}{\beta} $ is a pole for these both functions. 

However, if we consider $ \tilde{S}_{-\beta}$ as $ (-\beta)$-shift, the table above changes a little bit:  $ \tilde{S}_{-\beta}$ is coded if only if $ \beta $ is greater than or equal to the golden ratio $ \gamma_0$. But, if $ \beta $ less than the golden ratio, the systems $S_{-\beta} $ or $ \tilde{S}_{-\beta}$ contains coded sub-shift with maximal entropy $ \log \beta$. 



  \bibliographystyle{plain} 
 \bibliography{biblio(1)}

\end{document}